\documentclass{article}
\usepackage{geometry}
\usepackage[english]{babel}
\usepackage[utf8]{inputenc}
\usepackage{fontenc}
\usepackage{newlfont}
\usepackage{indentfirst}
\usepackage{fancyhdr}
\usepackage{textcomp}
\usepackage{graphicx}
\usepackage{amsmath}
\usepackage{amssymb}
\usepackage{amsthm}
\usepackage{mathrsfs}
\usepackage{float}
\usepackage[dvips]{hyperref}
\theoremstyle{plain}
\newtheorem{Teo}{Theorem}[section]
\newtheorem{Prop}[Teo]{Proposition} 
 
\newtheorem{Lem}[Teo]{Lemma}
\theoremstyle{definition}
\newtheorem{Def}[Teo]{Definition}

\newtheorem{Rem}[Teo]{Remark}

\title{Bellman function for Hardy's inequality over dyadic trees}
\author{Cavina Michelangelo}
\begin{document}

\maketitle
\begin{abstract}
In this article we use the Bellman function technique to characterize the measures for which the weighted Hardy's inequality holds on dyadic trees. We enunciate the (dual) Hardy's inequality over the dyadic tree and we use the associated ``Burkholder-type" function to find the main inequality required in the proof of the result. We enunciate a ``Bellman-type" function satisfying the required properties and we use the Bellman function technique to prove the main result of this article. We prove the optimality of the associated constant with an explicit example of an extremal family of maps. We also give an explicit interpretation of the corresponding Bellman function in terms of the theory of stochastic optimal control.
\end{abstract}
\tableofcontents
\section*{Introduction}
\addcontentsline{toc}{section}{Introduction}
The weighted Hardy's inequality on trees was initially studied for its applications to the theory of holomorphic function spaces, but it is an interesting topic on its own. The weighted dyadic Hardy's inequality was studied (see \cite{ARS1} and \cite{ARS2}) to characterize Carleson measures for analytic Besov space.\\
In this article we study the problem and solve it for the general case $1<p<+\infty$ and we prove that the inequality holds with constant $C(p)=\big(p/(p-1)\big)^p=\big(p'\big)^p$. See theorem \ref{Theorem in intro}. \\

But for the best constant, our characterization of the (dual) dyadic Hardy's inequality is not new, see \cite{ARS1}. The proof we give is new and it is inspired by the linear case given in \cite{AHMV}. The weighted dyadic Hardy's inequality can be characterized by other equivalent, but different conditions. For instance a capacitary characterization can be given, using the Maz'ya theory, see \cite{Mazya}.\\

In the past twenty years several results of this kind have been proved using the Bellman function method. The ideas behind the Bellman function technique were inspired by \cite{Burkholder}, see also \cite{Banuelos} . The expository article \cite{NTV} investigates the connection between the Bellman function technique in dyadic analysis and Bellman functions from the theory of stochastic control. The seminal article \cite{NT} presents a thorough exposition about the Bellman function technique and its applications. The article \cite{AHMV} solves the problem for the case $p=2$, and the Bellman function in our article is equal to the Bellman function used in \cite{AHMV} when we set $p=2$. The Bellman function we use in this article is very similar, but not completely, to the one used for the proof of the dyadic Carleson embedding theorem, see \cite{JL}, and for the Carleson embedding theorem the same constant $\big(p'\big)^p$ is sharp.\\

We are now going to present the results in this work.\\
Given the interval $I_0=[0,1]$ we denote by $\mathscr D(I_0)$ the standard dyadic tree structure of real intervals $I\subseteq I_0$.  We consider the maps 
\begin{align*}
I&\longmapsto \alpha_I \in \mathbb R^+,\\
I&\longmapsto \lambda_I \in \mathbb R^+,\\
I&\longmapsto\phi(I)\in \mathbb R^+,
\end{align*}
where we can read $\{\alpha_I\}$ as a choice of weights, $\{\lambda_I\}$ as a measure and $\{\phi(I)\}$ as a function over the dyadic tree.\\
The main result of this work is the following one:
\begin{Teo}\label{Theorem in intro}
Let $I_0$ be a real interval. Let $\{\alpha_I\}$ and $\{\lambda_I\}$ be a choice of weights and measure. If 
\begin{equation}\label{hypothesis intro}
\frac{1}{|I|}\sum_{K \subseteq I} \alpha_K\bigg(\frac{1}{|K|}\sum_{J \subseteq K} \lambda_J\bigg)^p \leq \frac{1}{|I|}\sum_{K \subseteq I}\lambda_I < +\infty \quad \forall I \in \mathscr D( I_0)
\end{equation}
is satisfied, then the dual weighted dyadic Hardy's inequality holds for $\{\alpha_I\}$ and $\{\lambda_I\}$, i.e.
\begin{equation}\label{thesis intro}
\frac{1}{|I_0|}\sum_{I \subseteq I_0} \alpha_I\bigg(\frac{1}{|I|}\sum_{J \subseteq I}\phi(J) \lambda_J^{\frac{1}{p'}}\bigg)^p \leq C(p)\frac{1}{|I_0|}\sum_{I \subseteq I_0}\phi(I)^p 
\end{equation}
for all $\phi \in l^p(\mathscr D(I_0))$. Here $C(p)$ is the constant
\[C(p)=\bigg(\frac{p}{p-1}\bigg)^p=\big(p'\big)^p.\]
The constant $C(p)$ in the inequality (\ref{thesis intro}) is sharp.\\
Moreover, if the inequality (\ref{thesis intro}) holds with constant $C(p)=1$ then the inequality (\ref{hypothesis intro}) holds for any $I \in \mathscr D(I_0)$ by choosing $\phi(I)=\lambda_I^{\frac{1}{p}}$ and by rescaling $I_0$ over $I$.
\end{Teo}

By setting $\eta(I)=\phi(I)\lambda_I^{-\frac{1}{p}}$ and $\omega_I^{1-p}=\frac{\alpha_I}{|I|^p}$ we rewrite the inequality (\ref{thesis intro}) in the form
\begin{equation*}
\sum_{I \subseteq I_0} \omega_I^{1-p}\bigg(\sum_{J \subseteq I}\eta(J)\lambda_J\bigg)^p \leq C(p)\sum_{I \subseteq I_0}\eta(I)^p \lambda_I,
\end{equation*}
which, by duality, is equivalent to the weighted dyadic Hardy's inequality
\begin{equation*}
\sum_{I \subseteq I_0} \lambda_I \bigg(\sum_{J \supseteq I}\psi(J)\bigg)^{p'} \leq C(p)\sum_{I \subseteq I_0}\psi(I)^{p'} \omega(I) \quad \text{for all } \psi  \in l^{p'}(\mathscr D(I_0)).
\end{equation*}

The function we use to prove our main result is
\[\mathcal B(F,f,A,v)=\bigg(\frac{p}{p-1}\bigg)^p F-\frac{p^p}{p-1}\frac{f^p}{(A+(p-1)v)^{p-1}}\]
defined over the domain
\[\mathcal D=\{(F,f,A,v) \in \mathbb R^4  \mid  F>0, f> 0, A>0,v>0,v\geq A, f^p \leq F v^{p-1}\}.\]
The properties of $\mathcal B$ are stated in subsection 2.2.\\

The function $\mathcal B$ can be interpreted as the solution of a Hamilton-Jacobi-Bellman equation associated to a stochastic problem of optimal control, which we will state in the article.\\

Our article is structured as follows.\\
In section 1 we show how to connect the dyadic setting to the analytical setting for Theorem \ref{Theorem in intro}. We enunciate the theorem in subsection 1.1, we define the ``Burkholder-type" function $B$ associated to Theorem \ref{Theorem in intro} in subsection 1.2, and we enunciate and prove the associated main inequality for $B$ in subsection 1.3.\\

In section 2 prove Theorem \ref{Theorem in intro} with the Bellman function method. We show that Theorem \ref{Theorem in intro} follows from the existence of a function of the ``Bellman-type" in the subsection 2.1, we show that the function $\mathcal B$ satisfies the required properties in the subsection 2.2, and we prove that the constant $C(p)$ is sharp for Theorem \ref{Theorem in intro} in the subsection 2.3.\\

In section 3 we state a stochastic optimal control problem whose solution is given by the Bellman function used throughout the paper. This gives a direct probabilistic interpretation to our function. We show a natural way to transition from a dyadic inequality to a Hamilton-Jacobi-Bellman inequality in the subsection 3.1, we define a stochastic optimal control problem whose solution is a Bellman function that satisfies the required Hamilton-Jacobi-Bellman inequality in the subsection 3.2, and we prove that the Bellman function associated to the stochastic optimal control problem we defined is equal to the function $\mathcal B$ in subsection 3.3.\\

\textbf{Keywords}--- Ahlfors-regular metric spaces, imbeddings of metric spaces.\\
\textbf{Mathematical Subject Classification 2020}--- 30H25, 30L05, 49L99.

\section{Hardy's inequality}
In this section we state the main result of this paper and we show a way to localize the problem in $\mathbb R^n$ using a ``Burkholder-type" function. We also find the main inequality associated to the problem, which will be used to prove the main result.
\subsection{Inequality over the dyadic tree}
Let $\mathscr D(I_0)$ be the dyadic tree over $I_0=[0,1]$.  Let $I \in \mathscr D(I_0)$, let
\begin{equation*}
\varphi:\mathscr D(I_0) \longrightarrow \mathbb R.
\end{equation*}
We denote the sum of the values $\varphi(J)$ for $J \in \mathscr D(I)$ by
\begin{equation*}
\sum_{J \in \mathscr D(I)} \varphi(J)=: \sum_{J \subseteq I} \varphi(J).
\end{equation*}
 Let $\Lambda$ be a positively valued measure over the dyadic tree defined as follows: for each node $I \in \mathscr D(I_0)$
\[\mathscr D(I_0) \ni I \longmapsto \lambda_I \in \mathbb R^+.\]
We define the following objects as follows:
\begin{align*}
\Lambda(I)&=\sum_{K \subseteq I}\lambda_K,\\
(\Lambda)_I &= \frac{1}{|I|} \sum_{K \subseteq I}\lambda_K=\frac{1}{|I|} \Lambda(I),\\
\int_I \phi \;d \Lambda &= \sum_{K \subseteq I} \phi(K) \lambda_K,\\
(\phi \Lambda)_I&=\frac{1}{|I|}\sum_{K \subseteq I} \phi(K)\lambda_K=\frac{1}{|I|}\int_I \phi \; d \Lambda.
\end{align*}
Our goal is to prove Theorem (1.3) in the article \cite{AHMV} in the general case $p \neq 2$.
\begin{Teo}[Dual weighted Hardy's inequality for dyadic trees]\label{hardy theorem}
Let $\mathscr D(I_0)$ be the dyadic tree originating at $I_0$, let $\{\alpha_I\}_{I \subseteq I_0}$ be a sequence of positive numbers. Let $\Lambda: \mathscr D(I_0) \rightarrow \mathbb R^+$ be a positive measure over the dyadic tree. Let $p$ be a real number such that $1<p<+\infty$. If the inequality
\begin{equation}\label{hypothesis}
\frac{1}{|I|}\sum_{K \subseteq I} \alpha_K(\Lambda)_K^p \leq (\Lambda)_I <+\infty \quad \forall I \in \mathscr D( I_0)
\end{equation}
is satisfied, then
\begin{equation}\label{thesis}
\frac{1}{|I_0|}\sum_{I \subseteq I_0} \alpha_I(\phi \Lambda^{\frac{1}{p'}})_I^p \leq C(p)(\phi^p)_{I_0},
\end{equation}
for all $\phi: \mathscr  D(I_0) \rightarrow \mathbb R^+$ such that $\phi \in l^p(\mathscr D(I_0))$.\\
Here $\frac{1}{p}+\frac{1}{p'}=1$, $C(p)=\big(p/(p-1)\big)^p=({p'})^p$ is a constant depending only on $p$, and
\begin{equation*}
(\phi \Lambda^{\frac{1}{p'}})_I=\frac{1}{|I|}\sum_{K\subseteq I} \phi(K) \lambda_K^{\frac{1}{p'}}, \quad (\Lambda^p)_{I_0}=\frac{1}{|I_0|}\sum_{I \subseteq I_0}\lambda_I^p.
\end{equation*}
\end{Teo}
We will prove this theorem using the Bellman function's method. 

\subsection{From dyadic inequality to real function}
In this subsection we define a new function using a tecnique developed by Burkholder, following the standard approach in the Bellman function's method.\\
Let $p \in \mathbb R$, $1<p<+\infty$. We consider the domain
\begin{equation}
\mathcal D:=\bigg \{ (F,f,A,v) \in \mathbb R^4 \mid  F>0, f>0, A > 0, v > 0, v \geq A,  f^p\leq F v^{p-1} \bigg\}.
\end{equation}
Let us consider maps
\begin{align*}
\alpha:I&\longmapsto \alpha_I \in \mathbb R^+,\\
\Lambda:I&\longmapsto \lambda_I \in \mathbb R^+,\\
\phi:I&\longmapsto\phi(I)\in \mathbb R^+,
\end{align*}
such that $\phi \in l^p(\mathscr D(I_0))$, and such that the hypothesis
\begin{equation}
\frac{1}{|I|}\sum_{K \subseteq I} \alpha_K(\Lambda)_K^p \leq (\Lambda)_I <+\infty \quad \forall I \in \mathscr D( I_0)
\end{equation}
holds. Let $I \in \mathscr D(I_0)$. Let us define the map
\begin{equation}\begin{aligned}
\Psi_I:(\alpha,\Lambda,\phi) \longmapsto (F_I,f_I,A_I,v_I) \in \mathbb R^4,
\end{aligned}\end{equation}
where
\begin{equation}\begin{aligned}
F_I=&\frac{1}{|I|}\sum_{J \subseteq I} \phi(J)^p,\\
f_I=&\frac{1}{|I|}\sum_{J \subseteq I} \phi(J)\lambda_J^{\frac{1}{p'}},\\
A_I=&\frac{1}{|I|} \sum_{J \subseteq I} \alpha_J \bigg(\frac{1}{|J|}\sum_{K \subseteq J} \lambda_K\bigg)^p, \\
v_I=&\frac{1}{|I|}\sum_{J \subseteq I} \lambda_J.
\end{aligned}\end{equation}
\begin{Def}[Burkholder-type function]
Let $(F,f,A,v) \in \mathbb R^4$. Let us consider the function
\begin{equation}\begin{aligned}\label{B function for Carleson Embedding}
B(F,f,A,v):=\frac{1}{|I|} \sup_{\underset{\Psi_I(\alpha,\Lambda,\phi)= (F,f,A,v)}{\alpha,\Lambda,\phi \text{ such that}}} \bigg\{ \sum_{J \subseteq I} \alpha_J\bigg(\frac{1}{|J|} \sum_{K \subseteq J}\phi(K) \lambda_K^{\frac{1}{p'}}\bigg)^p\bigg\},
\end{aligned}\end{equation}
where the quantity
\begin{equation*}
\sum_{J \subseteq I} \alpha_J\bigg(\frac{1}{|J|} \sum_{K \subseteq J}\phi(K) \lambda_K^{\frac{1}{p'}}\bigg)^p
\end{equation*}
is the quantity that appears in the left hand side of the thesis (\ref{thesis}) of Theorem \ref{hardy theorem}. We say that $B$ is the Burkholder-type function associated to Theorem \ref{hardy theorem}.\\
We observe that the function $B$ does not depend on the choice of the interval $I$, which is proved by rescaling the weights $\alpha$, the function $\phi$ and the measure $\Lambda$.
\end{Def}
To show that the function $B$ is well defined we prove the following proposition.
\begin{Prop}\label{map Psi is surjective}
For all $I \in \mathscr D(I_0)$ the map $\Psi_I$ is a surjective map which maps the set of $(\alpha,\Lambda,\phi)$ such that $\phi \in l^p(\mathscr D(I_0))$ and satisfying (\ref{hypothesis}) onto the domain $\mathcal D$.
\end{Prop}
\begin{proof}

Let $I \in \mathscr D(I_0)$. We show that $\Psi_I$ maps the set of 
\begin{align*}
\alpha:I&\longmapsto \alpha_I \in \mathbb R^+,\\
\Lambda:I&\longmapsto \lambda_I \in \mathbb R^+,\\
\phi:I&\longmapsto\phi(I)\in \mathbb R^+, \quad \phi \in l^p(\mathscr D(I_0))
\end{align*}
satisfying (\ref{hypothesis}) into the domain $\mathcal D$. By definition of $\Psi_I$ it follows that
\begin{equation}\label{first step Psi maps into}
\Psi_I(\alpha,\Lambda,\phi)=(F_I,f_I,A_I,v_I)\in \{(F,f,A,v) \in \mathbb R^4 \mid F>0,f>0,A>0,v>0\}.
\end{equation}
Moreover, by definition of $A_I$ and $v_I$ and by (\ref{hypothesis}) we have
\begin{equation}\label{second step Psi maps into}
A_I=\frac{1}{|I|}\sum_{K \subseteq I} \alpha_K\bigg(\frac{1}{|K|}\sum_{J \subseteq K} \lambda_J\bigg)^p \leq \frac{1}{|I|}\sum_{K \subseteq I}\lambda_K=v_I<+\infty,
\end{equation}
and by H\"older's inequality we have
\begin{align*}
f_I=& \frac{1}{|I|}\sum_{K \subseteq I}\phi(K) \lambda_K ^{\frac{1}{p'}}\leq\\
&\frac{1}{|I|^{\frac{1}{p}}}\bigg(\sum_{K \subseteq I}\phi(K)^p \bigg)^{\frac{1}{p}}   \frac{1}{|I|^{\frac{1}{p'}}}\bigg(\sum_{K \subseteq I}\lambda_K \bigg)^{\frac{1}{p'}}=\\
&(\phi^p)_I^{\frac{1}{p}}(\Lambda)_I^{\frac{1}{p'}}=F_I^{\frac{1}{p}} v_I ^{\frac{1}{p'}},
\end{align*}
which reformulates to
\begin{equation}\label{third step Psi maps into}
f_I^p\leq F_I \cdot v_I^{p-1}.
\end{equation}
So, from (\ref{first step Psi maps into}), (\ref{second step Psi maps into}) and (\ref{third step Psi maps into}) if follows that $\Psi_I(\alpha,\Lambda,\phi) \in \mathcal D$ for all $I \in \mathscr D(I_0)$ and for all $(\alpha, \Lambda,\phi)$, satisfying the required hypotheses.\\
The remaining part of the proof shows that $\Psi_I$ is a surjective map. The surjectivity of $\Psi_I$ is an interesting fact on its own, but it is not required for the proof of Theorem \ref{hardy theorem}, and it can be skipped. Let $I \in \mathscr D(I_0)$. Let $(F,f,A,v) \in \mathcal D$ be arbitrary. Up to rescaling $I$ on $I_0$ it is sufficient to show that there exist 
\begin{align*}
I&\longmapsto \alpha_I \in \mathbb R^+,\\
I&\longmapsto \lambda_I \in \mathbb R^+,\\
I&\longmapsto \phi(I)\in \mathbb R^+,
\end{align*}
such that the hypothesis
\begin{equation}\label{hypothesis in optimal domain proof}
\frac{1}{|I|}\sum_{K \subseteq I} \alpha_K\bigg(\frac{1}{|K|}\sum_{J \subseteq K} \lambda_J\bigg)^p \leq \frac{1}{|I|}\sum_{K \subseteq I}\lambda_I < +\infty \quad \forall I \in \mathscr D( I_0)
\end{equation}
holds, and such that
\begin{equation}\begin{aligned}
F=&\frac{1}{\left|I_0\right|}\underset{I \subseteq I_0}{\sum}  \phi(I)^p,\\
f=&\frac{1}{\left|I_0\right|}\underset{I \subseteq I_0}{\sum}  \phi(I)\lambda_I^{\frac{1}{p'}},\\
A=&\frac{1}{\left|I_0\right|} \underset{I \subseteq I_0}{\sum}  \alpha_I \bigg(\frac{1}{|I|}\underset{K \subseteq I}{\sum}   \lambda_K\bigg)^p, \\
v=&\frac{1}{\left|I_0\right|}\underset{I \subseteq I_0}{\sum}   \lambda_I.
\end{aligned}\end{equation}
Let us define the parameters
\begin{equation}\label{definition of parameters P1 and P2}
P_1:=\frac{ f}{ F^{\frac{1}{p}}v^{\frac{1}{p'}}}, \quad P_2:=\frac{A}{v}.
\end{equation}
By the definition of $\mathcal D$ we have $v\geq A$, $f^p \leq F v^{p-1}$, so we get
\begin{equation}
0<P_1\leq 1, \quad 0<P_2 \leq 1.
\end{equation}
Let $\phi_0>0$, $\lambda_0>0$,  $x_1 \in \left(0,\frac{1}{2}\right)$, $x_2 \in \left(0,\frac{1}{2^p}\right)$ to be chosen later. We define
\begin{equation}\begin{aligned}
\lambda_I:=&\lambda_0 \cdot  \left| I \right| ^{\log_{\frac 1 2}(x_1)},\\
\phi(I):=&\phi_0 \cdot \left| I \right| ^{\log_{\frac 1 2}(x_2)}.
\end{aligned}\end{equation}
The fact that $x_1 \in \left(0,\frac{1}{2}\right)$, $x_2 \in \left(0,\frac{1}{2^p}\right)$ yields the following computations:
\begin{equation}\begin{aligned}\label{expression for F in domain proof}
\frac{1}{\left|I_0\right|}\underset{I \subseteq I_0}{\sum}  \phi(I)^p =& \frac{1}{|I_0|}\underset{I \subseteq I_0}{\sum}  {\phi_0}^p \cdot \left| I \right| ^{p\log_{\frac 1 2}(x_2)}=\\
 &\frac{ {\phi_0}^p}{\left|I_0\right|}\sum_{k=0}^{+\infty} \left(\left|I_0\right|\frac{1}{2^k} \right)^{p\log_{\frac 1 2}(x_2)}\cdot 2^k=\\
 & {\phi_0}^p\left|I_0\right|^{p\log_{\frac 1 2}(x_2)-1}\sum_{k=0}^{+\infty} \left(\frac{1}{2^k} \right)^{\log_{\frac 1 2}({x_2}^p)}\cdot 2^k=\\
 & {\phi_0}^p\left|I_0\right|^{\log_{\frac 1 2}(2 {x_2}^p)}\sum_{k=0}^{+\infty} \left(2{x_2}^{p}\right)^k=\\
 & {\phi_0}^p\left|I_0\right|^{\log_{\frac 1 2}(2 {x_2}^p)} \frac{1}{1-2{x_2}^{p}}.
\end{aligned}\end{equation}
\begin{equation}\begin{aligned}\label{expression for v in domain proof}
\frac{1}{\left|I_0\right|}\underset{I \subseteq I_0}{\sum}   \lambda_I=&\frac{1}{\left|I_0\right|}\underset{I \subseteq I_0}{\sum}  \lambda_0 \cdot  \left| I \right| ^{\log_{\frac 1 2}(x_1)}=\\
 &\frac{\lambda_0}{\left|I_0\right|} \sum_{k=0}^{+\infty}  \left(\left| I_0 \right| \frac{1}{2^k}\right) ^{\log_{\frac 1 2}(x_1)}\cdot 2^k=\\
 & \lambda_0\left|I_0\right|^{\log_{\frac 1 2}(x_1)-1}\sum_{k=0}^{+\infty} \left(\frac{1}{2^k} \right)^{\log_{\frac 1 2}(x_1)}\cdot 2^k=\\
 & \lambda_0\left|I_0\right|^{\log_{\frac 1 2}(2 {x_1})}\sum_{k=0}^{+\infty} \left(2x_1\right)^k=\\
 & \lambda_0\left|I_0\right|^{\log_{\frac 1 2}(2 {x_1})} \frac{1}{1-2x_1}.
\end{aligned}\end{equation}
\begin{equation}\begin{aligned}\label{expression for f in domain proof}
\frac{1}{\left|I_0\right|}\underset{I \subseteq I_0}{\sum}  \phi(I)\lambda_I^{\frac{1}{p'}}=&\frac{1}{\left|I_0\right|}\underset{I \subseteq I_0}{\sum} \phi_0{\lambda_0}^{\frac{1}{p'}}\cdot  \left| I \right| ^{\frac{1}{p'}\log_{\frac 1 2}(x_1)+\log_{\frac 1 2}(x_2)}=\\
 &\frac{\phi_0{\lambda_0}^{\frac{1}{p'}}}{\left|I_0\right|} \sum_{k=0}^{+\infty}  \left(\left| I_0 \right| \frac{1}{2^k}\right)^{\frac{1}{p'}\log_{\frac 1 2}(x_1)+\log_{\frac 1 2}(x_2)} \cdot 2^k=\\
 &\phi_0{\lambda_0}^{\frac{1}{p'}}\left|I_0\right|^{\log_{\frac 1 2}\left({x_1}^{\frac{1}{p'}} x_2\right) -1} \sum_{k=0}^{+\infty}  \left(\frac{1}{2^k}\right)^{\log_{\frac 1 2}\left({x_1}^{\frac{1}{p'}} x_2\right)} \cdot 2^k=\\
 &\phi_0{\lambda_0}^{\frac{1}{p'}}\left|I_0\right|^{\log_{\frac 1 2}\left(2{x_1}^{\frac{1}{p'}} x_2\right) } \sum_{k=0}^{+\infty}\left({x_1}^{\frac{1}{p'}}x_2\right)^k \cdot 2^k=\\
 &\phi_0{\lambda_0}^{\frac{1}{p'}}\left|I_0\right|^{\log_{\frac 1 2}\left(2{x_1}^{\frac{1}{p'}} x_2\right) } \frac{1}{1-2{x_1}^{\frac{1}{p'}}x_2}.
\end{aligned}\end{equation}
So we get
\begin{equation}
0<\left|I_0\right|^{\log_{\frac 1 2}(2 {x_2}^p)} \frac{1}{1-2{x_2}^{p}}<+\infty, \quad 0<\left|I_0\right|^{\log_{\frac 1 2}(2 {x_1})} \frac{1}{1-2x_1}<+\infty.
\end{equation}
Now we choose
\begin{equation}\label{definition of phi zero}
\phi_0:=\left(\frac{F}{\left|I_0\right|^{\log_{\frac 1 2}(2 {x_2}^p)} \frac{1}{1-2{x_2}^{p}}}\right)^{\frac{1}{p}},
\end{equation}
\begin{equation}\label{definition of lambda zero}
\lambda_0:=\frac{v}{\left|I_0\right|^{\log_{\frac 1 2}(2 {x_1})} \frac{1}{1-2x_1}}.
\end{equation}
From (\ref{expression for F in domain proof}) and (\ref{definition of phi zero}) we get
\begin{equation}\label{F value proved}
\frac{1}{\left|I_0\right|}\underset{I \subseteq I_0}{\sum}  \phi(I)^p =F,
\end{equation}
from (\ref{expression for v in domain proof}) and (\ref{definition of lambda zero}) we get
\begin{equation}\label{v value proved}
\frac{1}{\left|I_0\right|}\underset{I \subseteq I_0}{\sum}   \lambda_I=v.
\end{equation}
Now we choose $x_1 \in \left(0,\frac{1}{2}\right)$, $x_2 \in \left(0,\frac{1}{2^p}\right)$ such that
\begin{equation*}
\frac{1}{\left|I_0\right|}\underset{I \subseteq I_0}{\sum}  \phi(I)\lambda_I^{\frac{1}{p'}}=f.
\end{equation*}
The last equation is equivalent to
\begin{equation*}
\frac{\frac{1}{\left|I_0\right|}\underset{I \subseteq I_0}{\sum}  \phi(I)\lambda_I^{\frac{1}{p'}}}{F^{\frac 1 p} v^{\frac{1}{p'}}}=\frac{f}{F^{\frac 1 p} v^{\frac{1}{p'}}}=P_1,
\end{equation*}
so, since we have (\ref{F value proved}) and (\ref{v value proved}), we are going find $x_1 \in \left(0,\frac{1}{2}\right)$, $x_2 \in \left(0,\frac{1}{2^p}\right)$ such that
\begin{equation}\label{P_1 equation to solve}
\frac{\frac{1}{\left|I_0\right|}\underset{I \subseteq I_0}{\sum}  \phi(I)\lambda_I^{\frac{1}{p'}}}{\left(\frac{1}{\left|I_0\right|}\underset{I \subseteq I_0}{\sum}  \phi(I)^p \right)^{\frac 1 p} \left(\frac{1}{\left|I_0\right|}\underset{I \subseteq I_0}{\sum}   \lambda_I\right)^{\frac{1}{p'}}}=P_1,
\end{equation}
for any arbitrary value $0<P_1\leq 1$.\\
We compute (\ref{P_1 equation to solve}) using (\ref{expression for F in domain proof}),  (\ref{expression for v in domain proof}) and  (\ref{expression for f in domain proof}) to get
\begin{equation}\begin{aligned}
P_1=&\frac{\phi_0{\lambda_0}^{\frac{1}{p'}}\left|I_0\right|^{\log_{\frac 1 2}\left(2{x_1}^{\frac{1}{p'}} x_2\right) } \frac{1}{1-2{x_1}^{\frac{1}{p'}}x_2}}{\left( {\phi_0}^p\left|I_0\right|^{\log_{\frac 1 2}(2 {x_2}^p)} \frac{1}{1-2{x_2}^{p}}\right)^{\frac{1}{p}}\left(\lambda_0\left|I_0\right|^{\log_{\frac 1 2}(2 {x_1})} \frac{1}{1-2x_1}\right)^{\frac{1}{p'}}}=\\
 &\frac{\phi_0{\lambda_0}^{\frac{1}{p'}}}{\phi_0{\lambda_0}^{\frac{1}{p'}}}\frac{\left|I_0\right|^{\log_{\frac 1 2}\left(2{x_1}^{\frac{1}{p'}} x_2\right) }}{\left|I_0\right|^{\frac{1}{p}\log_{\frac 1 2}(2 {x_2}^p)} \left|I_0\right|^{\frac{1}{p'}\log_{\frac 1 2}(2 {x_1})}}\frac{\left(1-2x_1\right)^{\frac{1}{p'}}\left(1-2{x_2}^{p}\right)^{\frac{1}{p}}}{1-2{x_1}^{\frac{1}{p'}}x_2}=\\
 &\frac{\left|I_0\right|^{\log_{\frac 1 2}\left(2{x_1}^{\frac{1}{p'}} x_2\right) }}{\left|I_0\right|^{\log_{\frac 1 2}\left(2^{\left(\frac{1}{p}+\frac{1}{p'}\right)}{x_1}^{\frac{1}{p'}}x_2\right)}}\frac{\left(1-2x_1\right)^{\frac{1}{p'}}\left(1-2{x_2}^{p}\right)^{\frac{1}{p}}}{1-2{x_1}^{\frac{1}{p'}}x_2}=\\
 & \frac{\left(1-2x_1\right)^{\frac{1}{p'}}\left(1-2{x_2}^{p}\right)^{\frac{1}{p}}}{1-2{x_1}^{\frac{1}{p'}}x_2}.
\end{aligned}\end{equation}
Let us define the function
\begin{equation}
g:\left(0,\frac 1 2\right) \times \left(0, \frac{1}{2^p}\right) \longrightarrow \mathbb R,
\end{equation}
\begin{equation}
g(x_1,x_2):=\frac{\left(1-2x_1\right)^{\frac{1}{p'}}\left(1-2{x_2}^{p}\right)^{\frac{1}{p}}}{1-2{x_1}^{\frac{1}{p'}}x_2}.
\end{equation}
To prove that for any arbitrary $0<P_1\leq 1$ there exist $x_1 \in \left(0,\frac{1}{2}\right)$, $x_2 \in \left(0,\frac{1}{2^p}\right)$ such that (\ref{P_1 equation to solve}) is satisfied, we prove that
\begin{equation}
(0,1] \subseteq \bigg \{ g(x_1,x_2) \mid x_1 \in \left(0,\frac{1}{2}\right), x_2 \in \left(0,\frac{1}{2^p}\right)\bigg\}.
\end{equation}
Let $x_1 \in \left(0,\frac{1}{2}\right)$, take $x_2={x_1}^{\frac{1}{p}}$. Then
\begin{equation}\label{value of g function in optimality of domain equal to 1}
g(x_1,x_2)=g(x_1,{x_1}^{\frac{1}{p}})=\frac{\left(1-2x_1\right)^{\frac{1}{p'}}\left(1-2 x_1\right)^{\frac{1}{p}}}{1-2{x_1}^{\frac{1}{p'}+\frac{1}{p}}}=1.
\end{equation}
To finish the proof we are going to show that
\begin{equation}\label{Liminf to prove in domain optimality}
\liminf_{(x_1,x_2) \rightarrow \left(\frac{1}{2},\frac{1}{2^p}\right)}g(x_1,x_2)=0.
\end{equation}
To prove this let us consider
\begin{equation}
\epsilon:=1-2x_1; \quad \delta:=1-2x_2^p,
\end{equation}
which gives us
\begin{equation}
x_1=x_1(\epsilon)=\frac{1-\epsilon}{2}; \quad x_2=x_2(\epsilon)=\left(\frac{1-\delta}{2}\right)^{\frac 1 p}.
\end{equation}
So we have $\epsilon \in (0,1)$, $\delta \in (0,1)$, and (\ref{Liminf to prove in domain optimality}) is equivalent to
\begin{equation}
\liminf_{(x_1,x_2) \rightarrow \left(\frac{1}{2},\frac{1}{2^p}\right)}g(x_1,x_2)=\liminf_{(\epsilon,\delta) \rightarrow \left(0,0\right)}\frac{\epsilon^{\frac{1}{p'}}\delta^{\frac{1}{p}}}{1-\left(1-\epsilon\right)^{\frac{1}{p'}}\left(1-\delta\right)^{\frac 1 p}}=0.
\end{equation}
Now, using little $o$ notation, by the Taylor formula we have
\begin{equation}
(1-x)^{\frac 1 p}=1-\frac x p +\underset{x\rightarrow 0}{o(x)}.
\end{equation}
So we get
\begin{equation}\begin{aligned}
g(x_1,x_2)=&\frac{\epsilon^{\frac{1}{p'}}\delta^{\frac{1}{p}}}{1-\left(1-\frac{\epsilon}{p'} +\underset{\epsilon\rightarrow 0}{o(\epsilon)}\right)\left(1-\frac \delta p +\underset{\delta\rightarrow 0}{o(\delta)}\right)}=\\
& \frac{\epsilon^{\frac{1}{p'}}\delta^{\frac{1}{p}}}{\frac{\epsilon}{p'}+\frac{\delta}{p}-\frac{\epsilon \delta}{pp'}-\underset{\epsilon\rightarrow 0}{o(\epsilon)}\left(1-\frac \delta p +\underset{\delta\rightarrow 0}{o(\delta)}\right)-\underset{\delta\rightarrow 0}{o(\delta)}\left(1-\frac{\epsilon}{p'} +\underset{\epsilon\rightarrow 0}{o(\epsilon)}\right)}.\\
\end{aligned}\end{equation}
Let $t>0$. Let us choose
\begin{equation}
\delta=\delta(\epsilon):=\epsilon^{1+t}.
\end{equation}
So we get
\begin{equation}\begin{aligned}
g(x_1,x_2)=\frac{\epsilon^{\frac{1}{p'}+\frac{1}{p}+\frac{t}{p}}}{\frac{\epsilon}{p'}+\frac{\epsilon^{1+t}}{p}-\frac{\epsilon^{1+\frac{1}{p}+\frac{1}{t}}}{pp'}-\underset{\epsilon\rightarrow 0}{o(\epsilon)}\left(1-\frac{\epsilon^{1+t}}{p} +\underset{\epsilon\rightarrow 0}{o(\epsilon^{1+t})}\right)-\underset{\epsilon\rightarrow 0}{o(\epsilon^{1+t})}\left(1-\frac{\epsilon}{p'} +\underset{\epsilon\rightarrow 0}{o(\epsilon)}\right)}.\\
\end{aligned}\end{equation}
Finally, we divide numerator and denominator by $\epsilon$, and we use little $o$ notation properties, to get
\begin{equation}\begin{aligned}
g(x_1,x_2)=\frac{\epsilon^{\frac{t}{p}}}{\frac{1}{p'}+\frac{\epsilon^{t}}{p}-\frac{\epsilon^{\frac{1}{p}+\frac{1}{t}}}{pp'}-\underset{\epsilon\rightarrow 0}{o(1)}\left(1-\frac{\epsilon^{1+t}}{p} +\underset{\epsilon\rightarrow 0}{o(\epsilon^{1+t})}\right)-\underset{\epsilon\rightarrow 0}{o(\epsilon^{t})}\left(1-\frac{\epsilon}{p'} +\underset{\epsilon\rightarrow 0}{o(\epsilon)}\right)}.\\
\end{aligned}\end{equation}

The denominator converges to $\frac{1}{p'}$ as $\epsilon \rightarrow 0$, while the numerator converges to 0 as $\epsilon \rightarrow 0$, which proves that, under the previous choices, we have
\begin{equation}
g(x_1(\epsilon),x_2(\epsilon)) \rightarrow 0 \quad \text{as } \epsilon \rightarrow 0.
\end{equation}
Since $g>0$ by construction, this entails 
\begin{equation}
\liminf_{(x_1,x_2) \rightarrow \left(\frac{1}{2},\frac{1}{2^p}\right)}g(x_1,x_2)=0.
\end{equation}
The function $g$ is a continuous function defined over the connected set $\{ (x_1,x_2) \mid x_1 \in \left(0,\frac{1}{2}\right), x_2 \in \left(0,\frac{1}{2^p}\right)\}$, so, since we also proved (\ref{value of g function in optimality of domain equal to 1}), it follows that
\begin{equation}
(0,1] \subseteq \bigg \{ g(x_1,x_2) \mid x_1 \in \left(0,\frac{1}{2}\right), x_2 \in \left(0,\frac{1}{2^p}\right)\bigg\},
\end{equation}
which is the equation we wanted to prove.\\
The last equation entails that (\ref{P_1 equation to solve}) is satisfied for any $0<P_1\leq 1$, so we proved that for all $F>0$, $f>0$, $v>0$ such that $f^p \leq F v^{p-1}$ then there exist $\lambda_0>0$, $\phi_0>0$, $0<x_1<\frac{1}{2}$, $0<x_2<\frac{1}{2^p}$ such that the maps
\begin{equation}\begin{aligned}
\lambda_I:=&\lambda_0 \cdot  \left| I \right| ^{\log_{\frac 1 2}(x_1)},\\
\phi(I):=&\phi_0 \cdot \left| I \right| ^{\log_{\frac 1 2}(x_2)}
\end{aligned}\end{equation}
satisfy the equations
\begin{equation}
\begin{aligned}
F=&\frac{1}{\left|I_0\right|}\underset{I \subseteq I_0}{\sum}  \phi(I)^p,\\
f=&\frac{1}{\left|I_0\right|}\underset{I \subseteq I_0}{\sum}  \phi(I)\lambda_I^{\frac{1}{p'}},\\
v=&\frac{1}{\left|I_0\right|}\underset{I \subseteq I_0}{\sum}   \lambda_I.
\end{aligned}
\end{equation}
Now it is only left to prove that there exist
\begin{equation}
\alpha:\mathscr D(I_0) \longrightarrow \mathbb R^+
\end{equation}
such that we have
\begin{equation}
A=\frac{1}{\left|I_0\right|} \underset{I \subseteq I_0}{\sum}  \alpha_I \bigg(\frac{1}{|I|}\underset{K \subseteq I}{\sum}   \lambda_K\bigg)^p. 
\end{equation}
Given the previous definitions of $\lambda_I$ and $\phi_I$, let us consider
\begin{equation}\begin{aligned}
\alpha:\mathscr D(I_0) \longrightarrow& \mathbb R^+\\
I \longmapsto& \alpha_I:= P_2 \cdot \lambda_I \bigg(\frac{1}{|I|}\underset{K \subseteq I}{\sum}   \lambda_K\bigg)^{-p},
\end{aligned}\end{equation}
where $P_2$ is the parameter defined in (\ref{definition of parameters P1 and P2}). This expression is well defined because, by construction of $\lambda_I$, we have 
\begin{equation}
\frac{1}{|I|}\underset{K \subseteq I}{\sum}   \lambda_K>0 \quad \forall I \in \mathscr D(I_0).
\end{equation}
From this definition, for all $I \in \mathscr D(I_0)$, we get
\begin{equation}
\frac{1}{\left|I\right|} \underset{J \subseteq I}{\sum}  \alpha_J \bigg(\frac{1}{|J|}\underset{K \subseteq J}{\sum}   \lambda_K\bigg)^p=P_2 \frac{1}{\left|I\right|} \underset{J \subseteq I}{\sum} \lambda_J\leq  \frac{1}{\left|I\right|} \underset{J \subseteq I}{\sum} \lambda_J<+\infty,
\end{equation}
which means that the hypothesis (\ref{hypothesis in optimal domain proof}) is satisfied. Moreover, for $I=I_0$, we get
\begin{equation}
\frac{1}{\left|I_0\right|} \underset{I \subseteq I_0}{\sum}  \alpha_I \bigg(\frac{1}{|I|}\underset{K \subseteq I}{\sum}   \lambda_K\bigg)^p=P_1 \frac{1}{\left|I_0\right|} \underset{I \subseteq I_0}{\sum} \lambda_I=P_1 \cdot v = A,
\end{equation}
which is the required inequality. So we proved that for all $(F,f,A,v) \in \mathcal D$ there exist 
\begin{align*}
I&\longmapsto \alpha_I \in \mathbb R^+,\\
I&\longmapsto \lambda_I \in \mathbb R^+,\\
I&\longmapsto \phi(I)\in \mathbb R^+,
\end{align*}
such that
\begin{equation*}
\frac{1}{|I|}\sum_{K \subseteq I} \alpha_K\bigg(\frac{1}{|K|}\sum_{J \subseteq K} \lambda_J\bigg)^p \leq \frac{1}{|I|}\sum_{K \subseteq I}\lambda_I < +\infty \quad \forall I \in \mathscr D( I_0)
\end{equation*}
holds, and such that
\begin{align*}
F=&\frac{1}{\left|I_0\right|}\underset{I \subseteq I_0}{\sum}  \phi(I)^p,\\
f=&\frac{1}{\left|I_0\right|}\underset{I \subseteq I_0}{\sum}  \phi(I)\lambda_I^{\frac{1}{p'}},\\
A=&\frac{1}{\left|I_0\right|} \underset{I \subseteq I_0}{\sum}  \alpha_I \bigg(\frac{1}{|I|}\underset{K \subseteq I}{\sum}   \lambda_K\bigg)^p, \\
v=&\frac{1}{\left|I_0\right|}\underset{I \subseteq I_0}{\sum}   \lambda_I,
\end{align*}
ending the proof.
\end{proof}
\begin{Rem}
Let us consider arbitrary maps
\begin{align*}
I&\longmapsto \alpha_I \in \mathbb R^+,\\
I&\longmapsto \lambda_I \in \mathbb R^+,\\
I&\longmapsto \phi(I)\in \mathbb R^+, 
\end{align*}
such that they satisfy the hypotheses of Theorem \ref{hardy theorem}. Let $I \in \mathscr D(I_0)$. Let $(F,f,A,v) \in \mathcal D$ be arbitrary. Then, the thesis  (\ref{thesis}) reformulates in the following way: for all $(F,f,A,v) \in \mathcal D$ we have
\begin{equation}
\frac{1}{|I_0|} \sum_{J \subseteq I_0} \alpha_J\bigg(\frac{1}{|J|} \sum_{K \subseteq J}\phi(K) \lambda_K^{\frac{1}{p'}}\bigg)^p\leq \left(p'\right)^p\frac{1}{|I_0|} \sum_{J \subseteq I_0}\phi_J^p.
\end{equation}
for all $(\alpha,\Lambda,\phi)$ such that they statisfy the hypotheses of Theorem \ref{hardy theorem} and such that $\Psi_{I_0}(\alpha,\Lambda,\phi)=(F,f,A,v)$.\\
However, from $\Psi_{I_0}(\alpha,\Lambda,\phi)=(F,f,A,v)$ it follows that
\begin{equation}
\frac{1}{|I_0|} \sum_{J \subseteq I_0}\phi_J^p = F,
\end{equation}
so, passing to the supremum, the thesis (\ref{thesis}) may be reformulated in the following way:
\begin{equation}
B(F,f,A,v)=\frac{1}{|I_0|} \sup_{\underset{\Psi_{I_0}(\alpha,\Lambda,\phi)= (F,f,A,v)}{\alpha,\Lambda,\phi \text{ such that}}} \bigg\{ \sum_{J \subseteq I_0} \alpha_J\bigg(\frac{1}{|J|} \sum_{K \subseteq J}\phi(K) \lambda_K^{\frac{1}{p'}}\bigg)^p\bigg\} \leq \left(p'\right)^p F
\end{equation}
for all $(F,f,A,v) \in \mathcal D$.
\end{Rem}
The previous observation yields the following proposition.
\begin{Prop}
The following statements are equivalent:
\begin{enumerate}
\item Theorem \ref{hardy theorem} holds.\\
\item The constant $C(p)=\left(p'\right)^p$ is the smallest constant for which we have
\begin{equation}
B(F,f,A,v) \leq C(p) \cdot F
\end{equation}
for all $(F,f,A,v) \in \mathcal D$.
\end{enumerate}
\end{Prop}
So, to prove Theorem \ref{hardy theorem}, we may instead study the properties of the function $B$. However, studying the properties of functions of the Burkholder-type like the function $B$ is, generally speaking, difficult. An idea from the work of Burkholder is to find a proper different function
\begin{equation*}
\mathcal B: \mathcal D \longmapsto \mathbb R
\end{equation*}
such that
\begin{enumerate}
\item $\mathcal B$ is defined over the same domain $\mathcal D$.\\
\item The function $\mathcal B$ satisfies the required bound, so in this case $\mathcal B$ satisfies
\begin{equation*}
0\leq \mathcal B(F,f,A,v) \leq \left(p'\right)^p F \quad \text{for all }(F,f,A,v) \in \mathcal D.
\end{equation*}
\item The function $\mathcal B$ satisfies an appropriate inequality, which is also satisfied by the function $B$, called main inequality.
\end{enumerate}
We will denote by function of the Bellman-type associated to Theorem \ref{hardy theorem} a function $\mathcal B$ satisfying the previous properties.\\
Given a function $\mathcal B$ with these properties, it is possible to prove that Theorem \ref{hardy theorem} holds, so our goal is to find the associated main inequality and a function $\mathcal B$ that satisfies the required properties.
\subsection{Main inequality}
The next lemma is about the main inequality, which will be the key to prove the dyadic Hardy's inequality.
\begin{Lem}\label{Main inequality lemma for B function}
For all $a,b,c\geq 0$, for all $(F,f,A,v),(F_-,f_-,A_-,v_-), (F_+,f_+,A_+,v_+) \in \mathcal D$ such that
\begin{align*}
&F=\tilde F+b^p ,  &f= \tilde f+ab,\\
&v=\tilde v+a^{p'},  &A =\tilde A+c,
\end{align*}
where
\begin{align*}
&\tilde F = \frac{1}{2}(F_-+F_+),  & \tilde f=\frac{1}{2}(f_-+f_+),\\
&\tilde v=\frac{1}{2}(v_-+v_+),  &\tilde A=\frac{1}{2}(A_-+A_+),
\end{align*}
 the function $B$ satisfies the main inequality
\begin{equation}
 B(F,f,A,v)-\frac{1}{2}\bigg[B(F_-,f_-,A_-,v_-)+ B(F_+,f_+,A_+,v_+)\bigg]\geq \frac{f^p}{v^p}c.
\end{equation}
\end{Lem}
\begin{proof}
 Given any fixed choice of $\tilde \alpha$, $\tilde \phi$, $\tilde \Lambda$ and $I\in \mathscr D(I_0)$, we consider the points
\begin{equation}\begin{aligned}
(F,f,A,v)=&\Psi_I(\tilde \alpha, \tilde \Lambda, \tilde \phi),\\
(F_+,f_+,A_+,v_+)=&\Psi_{I_+}(\tilde \alpha, \tilde \Lambda, \tilde \phi),\\
(F_-,f_-,A_-,v_-)=&\Psi_{I_-}(\tilde \alpha, \tilde \Lambda, \tilde \phi),
\end{aligned}\end{equation}
i.e.
\begin{equation}
\begin{cases}
F=\frac{1}{|I|}\underset{J \subseteq I}{\sum} \tilde \phi(J)^p,\\
f=\frac{1}{|I|}\underset{J \subseteq I}{\sum}  \tilde \phi(J)\tilde \lambda_J^{\frac{1}{p'}},\\
A=\frac{1}{|I|} \underset{J \subseteq I}{\sum}  \tilde \alpha_J \bigg(\frac{1}{|J|}\underset{K \subseteq J}{\sum}  \tilde \lambda_K\bigg)^p ,\\
v=\frac{1}{|I|}\underset{J \subseteq I}{\sum}  \tilde \lambda_J,
\end{cases}
\end{equation}
and $(F_{\pm},f_{\pm},A_{\pm},v_{\pm})$ are defined in the same way with $I_{\pm}$ instead of $I$.\\
We observe that these points satisfy the following equations:
\begin{align*}
F = \frac{1}{2}(F_-+F_+)+\frac{\tilde \phi(I)^p }{|I|},  \quad  &f=\frac{1}{2}(f_-+f_+)+\frac{\tilde \phi(I) \tilde \lambda_I^{\frac{1}{p'}}}{|I|}, \\
A=\frac{1}{2}(A_-+A_+)+\frac{\tilde \alpha_I \bigg(\frac{1}{|I|}\underset{J \subseteq I}{\sum}  \tilde \lambda_J\bigg)^p}{|I|}, \quad &v= \frac{1}{2}(v_-+v_+)+\frac{\tilde \lambda_I }{|I|}.
\end{align*} 
By setting
\begin{equation}\label{increments definition}
b=\frac{\tilde \phi(I) }{|I|^{\frac{1}{p}}},  \quad  c=\frac{\tilde \alpha_I \bigg(\frac{1}{|I|}\underset{J \subseteq I}{\sum}  \tilde \lambda_J\bigg)^p}{|I|}, \quad a=\frac{\tilde \lambda_I^{\frac{1}{p'}} }{|I|^{\frac{1}{p'}}},
\end{equation}
we get the following equations for the previous points:
\begin{equation}
\begin{matrix}\label{domain condition}
F = \frac{1}{2}(F_-+F_+)+b^p, \quad   &f=\frac{1}{2}(f_-+f_+)+ab, \\
\ \\
A=\frac{1}{2}(A_-+A_+)+c,\quad  &v= \frac{1}{2}(v_-+v_+)+a^{p'}.
\end{matrix}
\end{equation}
We compute
\begin{align*}
B(F,f,A,v)\geq\frac{1}{|I|} \sum_{J \subseteq I} \tilde \alpha_J\bigg(\frac{1}{|J|} \sum_{K \subseteq J}\tilde \phi(K) \tilde \lambda_K^{\frac{1}{p'}}\bigg)^p,
\end{align*}
which gives us
\begin{align*}
B(F,f,A,v)\geq&\frac{1}{|I|}\tilde \alpha_I\bigg(\frac{1}{|I|} \sum_{J \subseteq I}\tilde \phi(J) \tilde \lambda_J^{\frac{1}{p'}}\bigg)^p+\\
&\frac{1}{|I|} \sum_{J \in \mathscr D(I_+)} \tilde \alpha_J\bigg(\frac{1}{|J|} \sum_{K \subseteq J}\tilde \phi(K) \tilde \lambda_K^{\frac{1}{p'}}\bigg)^p+\\
&\frac{1}{|I|} \sum_{J \in \mathscr D(I_-)} \tilde \alpha_J\bigg(\frac{1}{|J|} \sum_{K \subseteq J}\tilde \phi(K) \tilde \lambda_K^{\frac{1}{p'}}\bigg)^p.
\end{align*}
Now we observe that
\begin{equation*}
\frac{\tilde \alpha_I}{|I|}=\frac{c}{v^p},
\end{equation*}
so the previous inequality becomes
\begin{equation}\begin{aligned}\label{step 2 in main inequality proof}
B(F,f,A,v)\geq\frac{f^p}{v^p}c+\frac{1}{|I|} \sum_{J \in \mathscr D(I_+)} \tilde \alpha_J\bigg(\frac{1}{|J|} \sum_{K \subseteq J}\tilde \phi(K) \tilde \lambda_K^{\frac{1}{p'}}\bigg)^p+\frac{1}{|I|} \sum_{J \in \mathscr D(I_-)} \tilde \alpha_J\bigg(\frac{1}{|J|} \sum_{K \subseteq J}\tilde \phi(K) \tilde \lambda_K^{\frac{1}{p'}}\bigg)^p.
\end{aligned}\end{equation}
Moreover, for any choice of  $\alpha$, $ \Lambda$ and  $\phi$ such that:
\begin{itemize}
\item $\Psi_I(\alpha,\Lambda,\phi)=(F,f,A,v)$,
\item $\Psi_{I_+}(\alpha,\Lambda,\phi)=(F_+,f_+,A_+,v_+)$,
\item  $\Psi_{I_-}(\alpha,\Lambda,\phi)=(F_-,f_-,A_-,v_-)$,
\end{itemize}
the following inequality holds:
\begin{equation}\begin{aligned}\label{step 3 in main inequality proof}
B(F,f,A,v)\geq\frac{f^p}{v^p}c+\frac{1}{|I|} \sum_{J \in \mathscr D(I_+)} \alpha_J\bigg(\frac{1}{|J|} \sum_{K \subseteq J} \phi(K) \lambda_K^{\frac{1}{p'}}\bigg)^p+\frac{1}{|I|} \sum_{J \in \mathscr D(I_-)}  \alpha_J\bigg(\frac{1}{|J|} \sum_{K \subseteq J} \phi(K)  \lambda_K^{\frac{1}{p'}}\bigg)^p.
\end{aligned}\end{equation}
So, by taking the supremum over all  $\alpha$,  $\Lambda$ and $ \phi$ for both the second and the third addend on the right hand side (using the fact that  $\alpha$,  $\Lambda$ and $ \phi$ can be ``independently" defined over $\mathscr D(I_+)$, $\mathscr D(I_-)$ and $I$) in the inequality (\ref{step 3 in main inequality proof}), we get
\begin{align*}
B(F,f,A,v)\geq&\frac{f^p}{v^p}c+\frac{1}{|I|}  \sup_{\underset{\Psi_{I_+}(\alpha,\Lambda,\phi)= (F_+,f_+,A_+,v_+)}{\alpha,\Lambda,\phi \text{ such that}}} \bigg\{ \sum_{J \subseteq I_+} \alpha_J\bigg(\frac{1}{|J|} \sum_{K \subseteq J}\phi(K) \lambda_K^{\frac{1}{p'}}\bigg)^p \bigg\}\\
&+\frac{1}{|I|}  \sup_{\underset{\Psi_{I_-}(\alpha,\Lambda,\phi)= (F_-,f_-,A_-,v_-)}{\alpha,\Lambda,\phi \text{ such that}}}\bigg\{ \sum_{J \subseteq I_-} \alpha_J\bigg(\frac{1}{|J|} \sum_{K \subseteq J}\phi(K) \lambda_K^{\frac{1}{p'}}\bigg)^p\bigg\}.
\end{align*}
Using the definition of the function $B$ and the fact that $|I|=2|I_+|=2|I_-|$ we get
\begin{equation*}
B(F,f,A,v)\geq \frac{f^p}{v^p} 
c+ \frac{1}{2}\bigg[ B(F_+,f_+,A_+,v_+)+B(F_-,f_-,A_-,v_-)\bigg].
\end{equation*}
The proof follows because the map $\Psi_I$ is surjective for all $I \in \mathscr D(I_0)$.
\end{proof}
\section{Bellman function's method}
In this section we prove Theorem \ref{hardy theorem} using the Bellman function's method. The proof will use the Bellman function
\begin{equation*}
\mathcal B(F,f,A,v)=\bigg(\frac{p}{p-1}\bigg)^p F-\frac{p^p}{p-1}\frac{f^p}{(A+(p-1)v)^{p-1}},
\end{equation*}
defined on the domain $\mathcal D$. We are going to show that the function $\mathcal B$ satisfies the man inequality from Lemma \ref{Main inequality lemma for B function}.\\
We are also going to show that the constant $\left(p'\right)^p$ is sharp for Theorem \ref{hardy theorem}. The proof will define a family of maps
\begin{align*}
\mathscr D(I_0) \ni I &\longmapsto \alpha_I(t) \equiv \alpha_I \in \mathbb R^+,\\
\mathscr D(I_0) \ni I &\longmapsto \lambda_I(t) \equiv \lambda_I \in \mathbb R^+,\\
\mathscr D(I_0) \ni I &\longmapsto \phi(t)(I) \equiv \phi(I) \in \mathbb R^+,
\end{align*}
for a proper family of indices $t \in T$, such that $\phi \in l^p(\mathscr D(I_0))$, such that
\begin{equation*}
\frac{1}{|I|} \sum_{K \subseteq I} \alpha_K(t)\left( \frac{1}{|K|} \sum_{J \subseteq K}\lambda_J(t)\right)^p \leq \frac{1}{|I|} \sum_{K \subseteq I} \lambda_K(t) < +\infty  \quad \forall I \in \mathscr D( I_0),
\end{equation*}
and such that for all $K(p)< \left(p'\right)^p$ there exist $\tilde t \in T$ such that
\begin{equation}
\frac{1}{|I_0|}\sum_{I \subseteq I_0}\alpha_I(\tilde t) \left(\frac{1}{|I|} \sum_{J \subseteq I}\phi_J(\tilde t) \lambda_J(\tilde t)^{\frac{1}{p'}}\right)^p>K(p) \cdot  \frac{1}{|I_0|} \sum_{J \subseteq I_0}\phi_J(\tilde t)^p.
\end{equation}
\subsection{The method to prove Hardy's inequality}
Now we will prove Theorem \ref{hardy theorem} using the Bellman function method. The argument is standard (see \cite{Burkholder}). We include the whole proof for the benefit of the reader.
\begin{Teo}[Bellman function's method for Theorem \ref{hardy theorem}]
Let $1<p<+\infty$, let $K(p)>0$ be an arbitrary constant depending only on $p$. Let $\mathcal D$ be the domain
\begin{equation*}
\mathcal D:=\bigg \{ (F,f,A,v) \in \mathbb R^4 \mid F>0, f>0, A > 0, v > 0, v \geq A,  f^p\leq F v^{p-1} \bigg\}.
\end{equation*}
Suppose that there exists a function
\begin{equation*}
\mathcal B:\mathcal D \longrightarrow \mathbb R
\end{equation*}
such that the following two properties hold:
\begin{enumerate}
\item $0<\mathcal B(F,f,A,v)\leq K(p) \cdot F \quad \text{for all }(F,f,A,v) \in \mathcal D,$\\
\item The main inequality holds for the function $\mathcal B$, i.e. for all $a,b,c\geq 0$, for all $(F,f,A,v)$, $(F_-,f_-,A_-,v_-)$, $ (F_+,f_+,A_+,v_+) \in \mathcal D$ such that
\begin{align*}
&F=\tilde F+b^p ,  &f= \tilde f+ab,\\
&v=\tilde v+a^{p'},  &A =\tilde A+c,
\end{align*}
where
\begin{align*}
&\tilde F = \frac{1}{2}(F_-+F_+),  & \tilde f=\frac{1}{2}(f_-+f_+),\\
&\tilde v=\frac{1}{2}(v_-+v_+),  &\tilde A=\frac{1}{2}(A_-+A_+),
\end{align*}
we have
\begin{equation}
 \mathcal B(F,f,A,v)-\frac{1}{2}\bigg[\mathcal B(F_-,f_-,A_-,v_-)+\mathcal B(F_+,f_+,A_+,v_+)\bigg]\geq \frac{f^p}{v^p}c.
\end{equation}
\end{enumerate}
Then, for all 
\begin{align*}
\alpha:I&\longmapsto \alpha_I \in \mathbb R^+,\\
\Lambda:I&\longmapsto \lambda_I \in \mathbb R^+,
\end{align*}
satisfying the hypothesis
\begin{equation}\label{hypothesis bellman function method theorem}
\frac{1}{|I|}\sum_{K \subseteq I} \alpha_K(\Lambda)_K^p \leq (\Lambda)_I <+\infty \quad \forall I \in \mathscr D( I_0),
\end{equation}
we have
\begin{equation}
\frac{1}{|I_0|}\sum_{I \subseteq I_0} \alpha_I (\phi \Lambda ^{\frac{1}{p'}})_{I}^p \leq K(p) \cdot (\phi^p)_{I_0}
\end{equation}
for all $\phi \in l^p(\mathscr D(I_0))$.
\end{Teo}
\begin{proof}
Let $\mathcal B$ be a function satisfying the required hypotheses. Let  $I \in \mathscr  D(I_0)$, we denote by $I_-\in \mathscr D(I_0)$ and $I_+\in \mathscr D(I_0)$ the two children of the node $I$.\\
Let 
\begin{align*}
\alpha:I&\longmapsto \alpha_I \in \mathbb R^+,\\
\Lambda:I&\longmapsto \lambda_I \in \mathbb R^+,
\end{align*}
such that they satisfy (\ref{hypothesis bellman function method theorem}), let $\phi \in l^p(\mathscr D(I_0))$.\\
For every $I \in \mathscr D(I_0)$ we consider
\begin{equation}
x_I=(F_I,f_I,A_I,v_I):=\Psi_I(\alpha,\Lambda,\phi),
\end{equation}
i.e.
\begin{align*}
v_I:=& (\Lambda)_I,\\
F_I:=& (\phi^p)_I,\\
f_I:=& (\phi \Lambda ^{\frac{1}{p'}})_I,\\
A_I:=& \frac{1}{|I|} \sum_{K \subseteq I} \alpha_K(\Lambda)^p_K.
\end{align*}
By Proposition \ref{map Psi is surjective} we have $(F_I,f_I,A_I,v_I) \in \mathcal D$ for all $I \in \mathscr D(I_0)$.\\
Now we define
\begin{align*}
a_I:=& \bigg( \frac{\lambda_I}{|I|} \bigg)^{\frac{1}{p'}},\\
b_I:=& \frac{\phi(I)}{|I|^{\frac{1}{p}}},\\
c_I:=& \frac{\alpha_I(\Lambda)^p_I}{|I|}.
\end{align*}
Let
\begin{align*}
x_{I_-}&=(F_{I_-},f_{I_-},A_{I_-},v_{I_-}):=\Psi_{I_-}(\alpha,\Lambda,\phi),\\
x_{I_+}&=(F_{I_+},f_{I_+},A_{I_+},v_{I_+}):=\Psi_{I_+}(\alpha,\Lambda,\phi).
\end{align*}
By computation, we get
\begin{align*}
v_I=&\frac{1}{|I|} \lambda_I + \frac{1}{2}(V_{I_-}+V_{I_+})=a_I^{p'} + \tilde v_I,\\
F_I=& \frac{1}{|I|} \phi(I)^p + \frac{1}{2} (F_{I_-}+F_{I_+})= b_I^p + \tilde F_I,\\
f_I=& \frac{\phi(I)\lambda_I^{\frac{1}{p'}}}{|I|}+\frac{1}{2}(f_{I_-}+f_{I_+})=a_I b_I+\tilde f_I,\\
A_I=& \frac{\alpha_I(\Lambda)^p_I}{|I|}+ \frac{1}{2}(A_{I_-}+A_{I_+})=c_I+\tilde A_I.
\end{align*}

So, for all choices of $\phi \in l^p(\mathscr D(I_0))$, $\alpha: \mathscr D(I_0) \rightarrow \mathbb R^+$, $\Lambda: \mathscr D(I_0) \rightarrow \mathbb R^+$, $I \in \mathscr D(I_0)$ satisfying the required hypotheses, the points
\[x_I := (F_I,f_I,A_I,v_I), \quad x_{I_-} := (F_{I_-},f_{I_-},A_{I_-},v_{I_-}), \quad x_{I_+} :=(F_{I_+},f_{I_+},A_{I_+},v_{I_+})\]
are eligible points for the main inequality satisfied by the function $\mathcal B$.\\
So we can apply the main inequality to get
\begin{align*}
|I|\frac{f_I^p}{v_I^p}c_I \leq& |I|\bigg[\mathcal B(x_I)-\frac{1}{2}\bigg(\mathcal B(x_{I_-})+\mathcal B(x_{I_+})\bigg)\bigg],\\
|I|\frac{f_I^p}{(\Lambda)_I^p}\frac{\alpha_I(\Lambda)_I^p}{|I|}\leq&|I|\mathcal  B(x_{I})-|I_-|\mathcal B(x_{I_-})-|I_+|\mathcal B(x_{I_+}),\\
\alpha_I f_I^p\leq&|I|\mathcal  B(x_{I})-|I_-|\mathcal B(x_{I_-})-|I_+|\mathcal B(x_{I_+}).
\end{align*}
Summing over all $I \in \mathscr D(I_0)$ and using the telescopic nature of the sum we get
\begin{equation}\label{Second last step proof bellman on tree}
\sum_{I \subseteq I_0} \alpha_I f_I^p \leq |I_0| \mathcal B(F_{I_0},f_{I_0},A_{I_0},v_{I_0}).
\end{equation}
However, by hypothesis we have $ \mathcal B(F,f,A,v)\leq K(p) \cdot F$ for all $(F,f,A,v) \in \mathcal D$, so we get
\begin{equation}
\sum_{I \subseteq I_0} \alpha_I f_I^p \leq |I_0| K(p) F_{I_0}.
\end{equation}
Now we recall that  $F_{I_0}=(\phi^p)_{I_0}$ and $f_{I}=(\phi \Lambda ^{\frac{1}{p'}})_{I}$, so we get
\begin{equation*}
\frac{1}{|I_0|}\sum_{I \subseteq I_0} \alpha_I (\phi \Lambda ^{\frac{1}{p'}})_{I}^p \leq K(p) \cdot (\phi^p)_{I_0},
\end{equation*}
which is the thesis, ending the proof.
\end{proof}

\subsection{The Bellman function for Hardy's inequality}
Let $p \in \mathbb R$, $1<p<+\infty$. We consider the function
\begin{equation}\label{bellman function}
\mathcal B(F,f,A,v)=\bigg(\frac{p}{p-1}\bigg)^p F-\frac{p^p}{p-1}\frac{f^p}{(A+(p-1)v)^{p-1}}
\end{equation}
defined over the domain
\[\mathcal D=\bigg \{ (F,f,A,v) \in \mathbb R^4 \mid F>0, f>0, A > 0, v > 0, v \geq A, f^p\leq F v^{p-1} \bigg\}.\]
Let $C(p)=\big(p/(p-1)\big)^p$. The function $\mathcal B$ has the following properties:
\begin{itemize}
\item[1)] $\mathcal B$ is a concave function defined over a convex domain.
\item[2)] $0\leq \mathcal B (F,f,A,v)\leq C(p) \cdot F $.
\end{itemize}
A proof of these properties can be found in the appendix.\\

Now we prove that the function $\mathcal B$ satisfies the main inequality associated to Theorem \ref{hardy theorem}.
\begin{Prop}\label{lemma}
 The function $\mathcal B$ satisfies
\begin{equation}\label{main inequality stronger}
\mathcal B(F,f,A,v)-\frac{1}{2}\bigg[\mathcal B(F_-,f_-,A_-,v_-)+\mathcal B(F_+,f_+,A_+,v_+)\bigg]\geq p^p\frac{f^p}{(A+v(p-1))^p}c,
\end{equation}
which, by using the fact that $v\geq A$, entails
\begin{equation}\label{main inequality 2}
\mathcal B(F,f,A,v)-\frac{1}{2}\bigg[\mathcal B(F_-,f_-,A_-,v_-)+\mathcal B(F_+,f_+,A_+,v_+)\bigg]\geq \frac{f^p}{v^p}c.
\end{equation}
where the inequality holds for all 
\begin{align*}
&F=\tilde F+b^p ,  &f= \tilde f+ab,\\
&v=\tilde v+a^{p'},  &A =\tilde A+c,\\
\end{align*}
and
\begin{align*}
&\tilde F = \frac{1}{2}(F_-+F_+),  & \tilde f=\frac{1}{2}(f_-+f_+),\\
&\tilde v=\frac{1}{2}(v_-+v_+),  &\tilde A=\frac{1}{2}(A_-+A_+),\\
\end{align*}
for every choice of $a\geq 0$, $b\geq 0$, $c\geq 0$. Here ${p'}$ is the real number such that $\frac{1}{p}+\frac{1}{p'}=1$.
\end{Prop}
\begin{proof}
We start by considering the telescopic sum
\begin{equation}\begin{aligned}\label{telescopic}\mathcal B(F,f,A,v)-\mathcal B(\tilde F,\tilde f,A-c,\tilde v)=&\mathcal B(F,f,A,v)-\mathcal B( F, f,A-c,v)+\\
 &\mathcal B( F,f,A-c, v)-\mathcal B( \tilde F, \tilde f, A-c,\tilde v).
\end{aligned}\end{equation}
Since the function $\mathcal B$ is concave and differentiable over a convex domain, we recall that a concave differentiable function's values are lower or equal to the values of any of its tangent hyperplanes. This entails that, for every $g$ concave and differentiable, for every choice of $x,x^* $ in the domain of the function $g$
\begin{equation}\label{originalconcavity}
g(x) - g(x^*) \leq \sum_{i=1}^4 \frac{\partial g(x^*)}{dx_i}(x_i-x^*_i).
\end{equation}
By changing the sign of (\ref{originalconcavity}) we get
\begin{equation}\label{concavity}
g(x^*) - g(x) \geq \sum_{i=1}^4 \frac{\partial g(x^*)}{dx_i}(x_i^*-x_i).
\end{equation}
So, when $g=\mathcal B$, $x=( F, f,  A, v)$, $x^*=(F,f,\tilde A,v)= ( F, f, A-c, v)$, the inequality  (\ref{concavity}) becomes 
\begin{equation}\label{increment}
\mathcal B(F,f,A,v)-\mathcal B(F,f,A-c,v)\geq p^p\frac{f^p}{(A+(p-1)v)^p}c.
\end{equation}
By combining (\ref{increment}) with (\ref{telescopic}) we get
\begin{equation}\label{telescopic2}
\mathcal B(F,f,A,v)-\mathcal B(\tilde F,\tilde f,A-c,\tilde v)\geq\mathcal B( F, f, A-c, v)-\mathcal B(\tilde F,\tilde f,A-c,\tilde v)+ (p-1)p^p\frac{f^p}{(A+(p-1)v)^p}c.
\end{equation}
Now we consider $g=\mathcal B$, $x=(\tilde F, \tilde f, A-c, \tilde v)$, $x^*=( F, f, A-c, v)$, so the inequality (\ref{concavity}) becomes
\begin{align*}
\mathcal B( F, f, A-c, v)- \mathcal B(\tilde F, \tilde f, A-c, \tilde v) \geq& \bigg(\frac{p}{p-1}\bigg)^p b^p-\frac{p^{p-1}}{p-1}\bigg(\frac{ f}{A-c+ (p-1)v}\bigg)^{p-1}ab+\\
 & (p-1) p^p\bigg(\frac{ f}{A-c+(p-1)v}\bigg)^{p}a^{p'}.
\end{align*}
Now let
\[y=\frac{ f}{A-c+ (p-1)v}.\]
We observe that $y> 0$ because $ f> 0$, $ v > 0$, $A-c > 0$ by definition of the domain of $\mathcal B$. So the last inequality can be rewritten in the form
\begin{equation*}
\mathcal B( F, f, A-c, v)- \mathcal B(\tilde F, \tilde f, A-c, \tilde v) \geq \bigg(\frac{p}{p-1}\bigg)^p b^p-\frac{p^{p+1}}{p-1}y^{p-1}ab+ (p-1)p^p y^{p}a^{p'}=:\phi(y).
\end{equation*}
Now we are going to prove that $\phi(y) \geq 0$ for all $y \geq 0$.\\
We observe that $\phi(y)=C(p)b^p \geq 0$ when $a=0$.\\
Now we assume $a>0$ and we compute the derivative of the function $\phi$:
\begin{align*}
\phi'(y)=p^{p+1} y^{p-2} \bigg((p-1)a^{p'} y -ab\bigg).
\end{align*}
So the derivative $\phi'(y)$ is such that $\phi'(y)\leq 0$ for $0 \leq y \leq \frac{b}{(p-1)a^{{p'}-1}}$, and  $\phi'(y)\geq 0$ for $y\geq \frac{b}{(p-1)a^{{p'}-1}}$, so $\tilde y=\frac{b}{(p-1)a^{{p'}-1}}$ is a point of absolute minimum for $\phi$, so as long as $\phi(\tilde y) \geq 0$ the inequality holds for all $y \geq 0$. So we compute
\begin{align*}
\phi(\tilde y)=& b^p-p\tilde y^{p-1}ab+ (p-1)\tilde y^{p}a^{p'}=\\
 & \bigg(\frac{p}{p-1}\bigg)^p b^p-\frac{p^{p+1}}{p-1}\bigg(\frac{b}{(p-1)a^{p'-1}}\bigg)^{p-1}ab+ (p-1)p^p\bigg(\frac{b}{(p-1)a^{p'-1}}\bigg)^{p}a^{p'}=\\
 & \bigg(\frac{p}{p-1}\bigg)^p b^p-\frac{p^{p+1}}{p-1}b^p\frac{1}{a^{pp'-p-p'}}+ \frac{p^p}{(p-1)^{p-1}} b^p\frac{1}{a^{pp'-p-p'}}.
\end{align*}
Now we recall that
\[\frac{1}{p}+\frac{1}{p'}=1; \quad p{p'}=p+{p'}.\]
So we get
\begin{align*}
\phi(\tilde y)=& \bigg(\frac{p}{p-1}\bigg)^p b^p-\frac{p^{p+1}}{p-1}b^p+ \frac{p^p}{(p-1)^{p-1}} b^p=\\
 &b^p \bigg(\frac{p}{p-1}\bigg)^p\bigg[1-p+p-1\bigg]=0.
\end{align*}
So the inequality $\phi(y)\geq 0$ holds for all $y \geq 0$, for every choice $a\geq 0$, $b \geq 0$, so the inequality (\ref{telescopic2}) becomes
\begin{equation}
\mathcal B(F,f,A,v)-\mathcal B(\tilde F,\tilde f,A-c,\tilde v)\geq p^p\frac{f^p}{(A+(p-1)v)^p}c.
\end{equation}
Now we observe that $(\tilde F,\tilde f, A-c,\tilde v)=(\tilde F,\tilde f,\tilde A,\tilde v)=\frac{1}{2}( (F_+,f_+,A_+,v_+)+ (F_-,f_-,A_-,v_-))$, so for the last step we use the fact that $\mathcal B$ is concave and we get
\[\mathcal B(F,f,A,v) - \frac{1}{2}\bigg[\mathcal B(F_+,f_+,A_+,v_+)+ \mathcal B(F_-,f_-,A_-,v_-)\bigg] \geq p^p\frac{f^p}{(A+(p-1)v)^p}c.\]
Finally, using the fact that $A\leq v$, we get the weaker version of the previous inequality
\[\mathcal B(F,f,A,v) - \frac{1}{2}\bigg[\mathcal B(F_+,f_+,A_+,v_+)+ \mathcal B(F_-,f_-,A_-,v_-)\bigg] \geq \frac{f^p}{v^p}c.\]
\end{proof}

\subsection{Sharpness of the constant}
In this subsection we prove that the constant
\[C(p):=\big(p'\big)^p=\bigg(\frac{p}{p-1}\bigg)^p\]
is sharp for Theorem \ref{hardy theorem}. The proof shows an example of an extremal family of maps
\begin{align*}
\alpha(t)\equiv \alpha:I&\longmapsto \alpha_I \in \mathbb R^+,\\
\Lambda(t)\equiv \Lambda:I&\longmapsto \lambda_I \in \mathbb R^+,\\
\phi(t)\equiv \phi:I&\longmapsto \phi(I)\in \mathbb R^+, \quad  \phi \in l^p(\mathscr D(I_0))
\end{align*}
for the constant $C(p)=\left(p'\right)^p$, i.e. a family of maps satisfying (\ref{hypothesis}), and such that
\begin{equation}\label{extremal family condition}
\lim_{t\rightarrow T}\frac{\frac{1}{|I_0|}\underset{I \subseteq I_0}{\sum}\alpha_I(t) \left(\frac{1}{|I|} \underset{J \subseteq I}{\sum}\phi_J(t) \lambda_J(t)^{\frac{1}{p'}}\right)^p}{\frac{1}{|I_0|} \underset{J \subseteq I_0}{\sum}\phi_J(t)^p}=\left(p'\right)^p.
\end{equation}
We spend few words on the idea behind the specific family of maps shown in the proof, which comes from the properties of the Burkholder-type function $B$ and the properties of our Bellman-type function $\mathcal B$. Generally speaking the Burkholder-type function associated to a problem is a solution to a proper Hamilton-Jacobi-Bellman equation, while any Bellman-type function associated to the same problem is a supersolution to the same Hamilton-Jacobi-Bellman equation (see Section 3), which yields that the Burkholder-type function $B$ is lower than or equal to any Bellman-type function associated to the same problem. So, if we consider an extremal family of maps $(\alpha(t),\Lambda(t),\phi(t))$, for $t_0<t<T$, and we set
\begin{equation}
(F(t),f(t),A(t),v(t)):= \Psi_{I_0}(\alpha(t),\Lambda(t),\phi(t)),
\end{equation}
from the idea that $\mathcal B \geq B$, and from $\mathcal B(F,f,A,v)\leq \left(p'\right)^p F$, we get 
\begin{align*}
\left(p'\right)^p=&\lim_{t\rightarrow T}\frac{\frac{1}{|I_0|}\underset{I \subseteq I_0}{\sum}\alpha_I(t) \left(\frac{1}{|I|} \underset{J \subseteq I}{\sum}\phi_J(t) \lambda_J(t)^{\frac{1}{p'}}\right)^p}{\frac{1}{|I_0|} \underset{J \subseteq I_0}{\sum}\phi_J(t)^p}\leq\\
&\lim_{t\rightarrow T}\left[\sup_{\underset{\Psi_{I_0}(\tilde \alpha,\tilde \Lambda,\tilde \phi)=(F(t),f(t),A(t),v(t))}{(\tilde \alpha,\tilde \Lambda,\tilde \phi) \text{ s.t.}}}\frac{\frac{1}{|I_0|}\underset{I \subseteq I_0}{\sum}{\tilde \alpha}_I \left(\frac{1}{|I|} \underset{J \subseteq I}{\sum}{\tilde \phi}_J {\tilde \lambda}_J^{\frac{1}{p'}}\right)^p}{F(t)}\right]=\\
&\lim_{t\rightarrow T}\frac{B(F(t),f(t),A(t),v(t))}{F(t)}\leq\\
&\lim_{t\rightarrow T}\frac{ \mathcal B(F(t),f(t),A(t),v(t))}{F(t)}\leq \\
&\left(p'\right)^p.
\end{align*}
Hence we looked for a proper family of maps $(\alpha(t),\Lambda(t),\phi(t))$ such that 
\begin{equation}
\lim_{t\rightarrow T}\frac{ \mathcal B(F(t),f(t),A(t),v(t))}{F(t)}= \left(p'\right)^p
\end{equation}
and we found an extremal family such that (\ref{extremal family condition}) holds.
\begin{Teo}
The constant $\big(p'\big)^p$ is sharp for Theorem \ref{hardy theorem}.
\end{Teo}
\begin{proof}
Let us consider $t_0<t<\frac 1 2$, for $0<t_0< \frac 1 2$, $t_0$ sufficiently close to $\frac 1 2$. Consider a function
\begin{equation}
g:\left(t_0,\frac 1 2 \right) \longmapsto \left(0,\frac{1}{2^p}\right)
\end{equation}
such that, using little $o$ notation, we have
\begin{equation}\label{little o condition for g optimal weights}
g\left(\frac 1 2 -t\right)=\frac{1}{2^p}+ \underset{t\rightarrow \left(\frac 1 2\right)^-}{o\left(\left|\frac 1 2 -t \right|\right)}.
\end{equation}
An example of a function $g$ with such properties is
\begin{equation}
g(t):= \frac{1}{2^p}-\left(\frac 1 2 -t\right)^s,
\end{equation}
for any fixed $s>1$, for $t_0=t_0(s)$ sufficiently close to $\frac 1 2$.\\
Let $\lambda_0>0$, $\phi_0>0$. Let us consider, for $t_0<t<\frac 1 2 $, the following maps
\begin{equation}\begin{aligned}
\lambda_I:=& \lambda_0 \cdot | I |^{\log_{\frac 1 2 }(t)}\quad \text{for }I \in \mathscr D(I_0),\\
\phi(I):=& \phi_0\cdot | I |^{\log_{\frac 1 2 }(g(t))} \quad \text{for }I \in \mathscr D(I_0) .
\end{aligned}\end{equation}
We observe that, by construction we have $0<2t<1$, so we compute
\begin{equation}\begin{aligned}
\frac{1}{|I|} \sum_{J \subseteq I}\lambda_J=&\frac{1}{|I|}  \sum_{J \subseteq I}\lambda_0 |J|^{\log_{\frac 1 2 }(t)}\\
=&\frac{1}{|I|} \sum_{k=0}^{+\infty} \lambda_0 \left(|I|\frac{1}{2^k}\right)^{\log_{\frac 1 2 }(t)}\cdot 2^k\\
=&\frac{1}{|I|} \lambda_0 |I|^{\log_{\frac 1 2 }(t)}\sum_{k=0}^{+\infty}t^k 2^k=\\
=& \lambda_0 |I|^{\log_{\frac 1 2 }(2t)}\frac{1}{1-2t},
\end{aligned}\end{equation}
which proves that
\begin{equation}\label{lambda in L^1 for optimal weights}
\frac{1}{|I|} \sum_{J \subseteq I}\lambda_J<+\infty.
\end{equation}
Moreover, by construction we have $0<2g(t)^p<1$, so we compute
\begin{equation}\begin{aligned}\label{l^p norm of phi in optimal weights}
\frac{1}{|I|} \sum_{J \subseteq I}\phi_J^p=&\frac{1}{|I|}  \sum_{J \subseteq I}\phi_0^p |J|^{p\log_{\frac 1 2 }(g(t))}\\
=&\frac{1}{|I|} \sum_{k=0}^{+\infty} \phi_0^p \left(|I|\frac{1}{2^k}\right)^{\log_{\frac 1 2 }(g(t)^p)}\cdot 2^k\\
=&\frac{1}{|I|} \phi_0^p |I|^{\log_{\frac 1 2 }(g(t)^p)}\sum_{k=0}^{+\infty}\left(g(t)^p\right)^k 2^k=\\
=& \phi_0^p |I|^{\log_{\frac 1 2 }(2g(t)^p)} \frac{1}{1-2g(t)^p},
\end{aligned}\end{equation}
which, taking $I=I_0$, proves that
\begin{equation}
\sum_{J \in \mathscr D(I_0)}\phi_J^p<+\infty.
\end{equation}
Now we define
\begin{equation}
\alpha:\mathscr D(I_0) \longrightarrow \mathbb R^+,
\end{equation}
such that
\begin{equation}
\alpha_K\left( \frac{1}{|K|} \sum_{J \subseteq K}\lambda_J\right)^p=\lambda_K\quad \text{for }K \in \mathscr D(I_0).
\end{equation}
So, by computation, we have 
\begin{equation}\begin{aligned}
\alpha_K:=&\lambda_K \left( \frac{1}{|K|} \sum_{J \subseteq K}\lambda_J\right)^{-p}\\
=& \lambda_0^{1-p} |K|^{\left(\log_{\frac 1 2 }(t)-p\log_{\frac 1 2 }(2t)\right)}(1-2t)^p\\
=& \lambda_0^{1-p} |K|^{\log_{\frac 1 2 }\left(\frac{1}{2^p} t^{(1-p)}\right)}(1-2t)^p.
\end{aligned}\end{equation}
From the definition of $\alpha$ and from (\ref{lambda in L^1 for optimal weights}) we get
\begin{equation}
\frac{1}{|I|} \sum_{K \subseteq I} \alpha_K\left( \frac{1}{|K|} \sum_{J \subseteq K}\lambda_J\right)^p \leq \frac{1}{|I|} \sum_{K \subseteq I} \lambda_K < +\infty,
\end{equation}
which is the required hypothesis. So we proved that $\alpha$, $\Lambda$ and $\phi$ satisfy the required hypotheses by Theorem \ref{hardy theorem}.\\
Now, since $0<2t^{\frac{1}{p'}}g(t)<1$, we may compute
\begin{equation}\begin{aligned}
\frac{1}{|I|} \sum_{J \subseteq I}\phi_J \lambda_J^{\frac{1}{p'}}=&\frac{1}{|I|}  \sum_{J \subseteq I}\phi_0 \lambda_0^{\frac{1}{p'}} |J|^{\log_{\frac 1 2 }(g(t))+\frac{1}{p'}\log_{\frac 1 2 }(t)}\\
=&\frac{1}{|I|}  \sum_{k=0}^{+\infty}\phi_0 \lambda_0^{\frac{1}{p'}} \left(|I|\frac{1}{2^k}\right)^{\log_{\frac 1 2 }(g(t))+\frac{1}{p'}\log_{\frac 1 2 }(t)}\cdot 2^k\\
=&\frac{1}{|I|} \phi_0 \lambda_0^{\frac{1}{p'}} |I|^{\log_{\frac 1 2 }\left(t^{\frac{1}{p'}}g(t)\right)} \sum_{k=0}^{+\infty} \left(2t^{\frac{1}{p'}}g(t)\right)^k\\
=& \phi_0 \lambda_0^{\frac{1}{p'}} |I|^{\log_{\frac 1 2 }\left(2t^{\frac{1}{p'}}g(t)\right)}\frac{1}{1-2t^{\frac{1}{p'}}g(t)}.\\
\end{aligned}\end{equation}
Now we compute the left hand side of the thesis:
\begin{equation}\begin{aligned}\label{left hand side thesis optimal weights}
 &\frac{1}{|I_0|}\sum_{I \subseteq I_0}\alpha_I \left(\frac{1}{|I|} \sum_{J \subseteq I}\phi_J \lambda_J^{\frac{1}{p'}}\right)^p=\\
 &\frac{1}{|I_0|}\sum_{I \subseteq I_0}\alpha_I\left(\phi_0 \lambda_0^{\frac{1}{p'}} |I|^{\log_{\frac 1 2 }\left(2t^{\frac{1}{p'}}g(t)\right)}\frac{1}{1-2t^{\frac{1}{p'}}g(t)}\right)^p=\\
 &\frac{1}{|I_0|}\phi_0^p \lambda_0^{p-1}\frac{1}{\left(1-2t^{\frac{1}{p'}}g(t)\right)^p}\sum_{I \subseteq I_0}\alpha_I |I|^{\log_{\frac 1 2 }\left(2^pt^{p-1}g(t)^p\right)}=\\
 &\frac{1}{|I_0|}\phi_0^p \lambda_0^{p-1}\frac{(1-2t)^p}{\left(1-2t^{\frac{1}{p'}}g(t)\right)^p}\sum_{I \subseteq I_0}\left(\lambda_0^{1-p} |I|^{\log_{\frac 1 2 }\left(\frac{1}{2^p} t^{(1-p)}\right)}\right)\left( |I|^{\log_{\frac 1 2 }\left(2^pt^{p-1}g(t)^p\right)}\right)=\\
&\frac{1}{|I_0|}\phi_0^p \frac{(1-2t)^p}{\left(1-2t^{\frac{1}{p'}}g(t)\right)^p}\sum_{I \subseteq I_0} |I|^{\log_{\frac 1 2 }\left(g(t)^p\right)}=\\
&\frac{1}{|I_0|}\phi_0^p \frac{(1-2t)^p}{\left(1-2t^{\frac{1}{p'}}g(t)\right)^p}\sum_{k=0}^{+\infty}\left(|I_0|\frac{1}{2^k}\right)^{\log_{\frac 1 2 }(g(t)^p)}\cdot 2^k=\\
&\phi_0^p  |I_0|^{\log_{\frac 1 2 }(2g(t)^p)} \frac{1}{1-2g(t)^p}\frac{(1-2t)^p}{\left(1-2t^{\frac{1}{p'}}g(t)\right)^p}. 
\end{aligned}\end{equation}
To finish the proof, we are going to show that
\begin{equation}
\lim_{t \rightarrow \left(\frac{1}{2}\right)^-}\frac{1}{|I_0|}\sum_{I \subseteq I_0}\alpha_I \left(\frac{1}{|I|} \sum_{J \subseteq I}\phi_J \lambda_J^{\frac{1}{p'}}\right)^p \left( \frac{1}{|I_0|} \sum_{J \subseteq I_0}\phi_J^p\right)^{-1} = \left(p'\right)^p.
\end{equation}
By computation, using (\ref{l^p norm of phi in optimal weights}) and (\ref{left hand side thesis optimal weights}), we get
\begin{equation}\begin{aligned}\label{object that converges to the optimal constant}
&\frac{1}{|I_0|}\sum_{I \subseteq I_0}\alpha_I \left(\frac{1}{|I|} \sum_{J \subseteq I}\phi_J \lambda_J^{\frac{1}{p'}}\right)^p \left( \frac{1}{|I_0|} \sum_{J \subseteq I_0}\phi_J^p\right)^{-1}=\\
&\phi_0^p  |I_0|^{\log_{\frac 1 2 }(2g(t)^p)} \frac{1}{1-2g(t)^p}\frac{(1-2t)^p}{\left(1-2t^{\frac{1}{p'}}g(t)\right)^p} \left( \phi_0^p |I|^{\log_{\frac 1 2 }(2g(t)^p)} \frac{1}{1-2g(t)^p}\right)^{-1}=\\
&\left(\frac{1-2t}{1-2t^{\frac{1}{p'}}g(t)}\right)^p.
\end{aligned}\end{equation}
So, to finish the proof, we are going to show that
\begin{equation}\label{limit to calculate optimal weights}
\lim_{t \rightarrow \left(\frac{1}{2}\right)^-} \frac{1-2t}{1-2t^{\frac{1}{p'}}g(t)}=p'.
\end{equation}
We consider the change of variable
\begin{equation}
x=1-2t; \quad t=\frac{1-x}{2}.
\end{equation}
By change of variable we get the equation
\begin{equation}
\lim_{t \rightarrow \left(\frac{1}{2}\right)^-} \frac{1-2t}{1-2t^{\frac{1}{p'}}g(t)}=\lim_{x \rightarrow 0^+} \frac{x}{1-2^{\frac{1}{p}}(1-x)^{\frac{1}{p'}}g\left(\frac{1-x}{2}\right)}.
\end{equation}
By the Taylor formula we have
\begin{equation}
(1-x)^{\frac{1}{p'}}=1-\frac{x}{p'}+\underset{x\rightarrow0^+}{o\left(\left|x \right|\right)},
\end{equation}
so we get
\begin{equation}\label{intermediate step limit optimal weights}
\lim_{t \rightarrow \left(\frac{1}{2}\right)^-} \frac{1-2t}{1-2t^{\frac{1}{p'}}g(t)}=\lim_{x \rightarrow 0^+} \frac{x}{1-2^{\frac{1}{p}}\left(1-\frac{x}{p'}+\underset{x\rightarrow0^+}{o\left(\left|x \right|\right)}\right)g\left(\frac{1-x}{2}\right)}.
\end{equation}
By change of variable, from (\ref{little o condition for g optimal weights}) we get 
\begin{equation}
g\left(\frac{1-x}{2}\right)=\frac{1}{2^p}+ \underset{x\rightarrow0^+}{o\left(\left|x \right|\right)},
\end{equation}
so equation (\ref{intermediate step limit optimal weights}) becomes
\begin{equation}\begin{aligned}\label{limit that gives the optimal constant}
\lim_{t \rightarrow \left(\frac{1}{2}\right)^-} \frac{1-2t}{1-2t^{\frac{1}{p'}}g(t)}=&\lim_{x \rightarrow 0^+} \frac{x}{1-2^{\frac{1}{p}}\left(1-\frac{x}{p'}+\underset{x\rightarrow0^+}{o\left(\left|x \right|\right)}\right)\left(\frac{1}{2^p}+ \underset{x\rightarrow0^+}{o\left(\left|x \right|\right)}\right)}\\
=&\lim_{x \rightarrow 0^+} \frac{x}{\frac{x}{p'}+\underset{x\rightarrow0^+}{o\left(\left|x \right|\right)}-2^{\frac{1}{p}} \left(1-\frac{x}{p'}+\underset{x\rightarrow0^+}{o\left(\left|x \right|\right)}\right)\underset{x\rightarrow0^+}{o\left(\left|x \right|\right)}}\\
=&\lim_{x \rightarrow 0^+} \frac{1}{\frac{1}{p'}+\underset{x\rightarrow0^+}{o\left(1\right)}-2^{\frac{1}{p}}\left(1-\frac{x}{p'}+\underset{x\rightarrow0^+}{o\left(\left|x \right|\right)}\right)\underset{x\rightarrow0^+}{o\left(1 \right)}}\\
=&p'.
\end{aligned}\end{equation}
So (\ref{object that converges to the optimal constant}) and (\ref{limit that gives the optimal constant}) yield that, for the given choice of $\alpha$, $\Lambda$ and $\phi$, we have
\begin{equation}
\lim_{t \rightarrow \left(\frac{1}{2}\right)^-}\frac{1}{|I_0|}\sum_{I \subseteq I_0}\alpha_I \left(\frac{1}{|I|} \sum_{J \subseteq I}\phi_J \lambda_J^{\frac{1}{p'}}\right)^p=\left(p'\right)^p  \frac{1}{|I_0|} \sum_{J \subseteq I_0}\phi_J^p,
\end{equation}
which proves that the constant $\left(p'\right)^p$ is sharp for Theorem \ref{hardy theorem}.
\end{proof}

\section{Stochastic approach to the problem}
We will now analyze this problem from the point of view of the theory of stochastic optimal control, and we will show that the function $\mathcal B$ can be interpreted as the Bellman function associated to a stochastic optimal control problem naturally related to the dyadic problem. In this section we use the same notations used in \cite{Oksendal}, chapter 11. See \cite{Krylov} for more details about the topic.\\

We are going to show the following theorem.
\begin{Teo}\label{dyadic bellman is stochastic theorem}
For all $x \in \mathcal D$ we have
\begin{equation}
\mathcal B(x)=g(x).
\end{equation}
\end{Teo}
Here $g$ is the Bellman function solution to the following stochastic optimal control problem associated to the inequality (\ref{main inequality stronger}).\\
Consider $x \in \mathcal D$, $u=(u_1,u_2,\dots,u_5) \in \mathbb R^5$ such that $u_5\geq0$. Let us define the payoff density
\[{\eta}^u(x):=p^p\bigg(\frac{x_2}{x_3+(p-1)x_4}\bigg)^p u_5.\]
Let $x \in \overline{\mathcal D}$. We define the bequest function
\[K(x)=\liminf_{\underset{y \in \mathcal D}{y\rightarrow x}}\mathcal B(y).\]
We remark that, for the definition of the stochastic Bellman function, we only need to define the bequest function $K$ on the boundary of the domain $\mathcal D$, however we follow the definition used in \cite{Oksendal}.\\
Let us define the coefficients
\[b(u,x):=(0,0,-u_5,0),\]
\[\sigma (u,x):=(u_1,u_2,u_3,u_4). \]
Let $\{u_t\}_{t\geq0}$ be a control such that $u_t(\omega) \in \{u \in \mathbb R^5 \; | \; u_5 \geq 0\}$. We consider the stochastic process $\{X_t\}=\{(F_t,f_t,A_t,v_t)\}$ solution to the following stochastic differential equation
\begin{equation}\label{process bellman def}
X_t=x_0+\int_0^t b(u_s,X_s)ds + \int_0^t \sigma(u_s,X_s)dB_s,
\end{equation}
where $x_0 \in \mathcal D$ is the  starting point, $\{B_t\}_{t\geq0}$ is a 1-dimensional Brownian motion and the domain of values of $X_t$ is the set $\mathcal D$. Let $\tau_{\mathcal D}$ be the first exit time for $\{X_t\}_{t\geq0}$ from $\mathcal D$, i.e.
\[\tau_{\mathcal D}(\omega):= \begin{cases}
 \inf \{ s>0 \mid X_s(\omega) \not \in \mathcal D\} \quad \text{if }  \{ s>0 \mid X_s(\omega) \not \in \mathcal D\} \neq \emptyset, \\
+\infty \quad \text{otherwise.}\end{cases}\]
The Bellman function associated to the problem is
\[g(x)=\sup_{\{u_t\}} E^x\bigg[\int_0^{\tau_{\mathcal D}} p^p\bigg(\frac{f_s}{A_s+(p-1)v_s}\bigg)^p u_5 \;ds +K(X_{\tau_{\mathcal D}})\chi_{\{\tau_{\mathcal D}<+\infty\}}\bigg],\]
where the supremum is taken over the set of controls $\{u_t\}_{t\geq0}$ satisfying proper measurability conditions and whose values range in the set $\{(u_1,u_2,u_3,u_4,u_5) \in \mathbb R^5\; | \; u_5 \geq 0\}$.\\

We observe that, for this result, we used the stronger version of the main inequality (\ref{main inequality stronger}) instead of the weaker version (\ref{main inequality 2}). By using a stronger main inequality we still get a Bellman function that can be used in the proof of Theorem \ref{hardy theorem} with the Bellman function method, however finding the solution to the problem associated to the weaker inequality (\ref{main inequality 2}) would require more work.\\

Now in the following subsections we show how we got to the stochastic optimal control problem and how we solved it.
\subsection{From the dyadic to the stochastic problem}
In this subsection we will show that the main inequality satisfied by the function $\mathcal B$ can be used to prove that $\mathcal B$ satisfies a differential inequality that will be the starting point from which we enunciate the stochastic optimal control problem having $\mathcal B$ as a solution.\\
We are going to recall the problem we are considering. Let $p \in \mathbb R$, $1<p<+\infty$. We consider the function
\begin{equation*}
\mathcal B(F,f,A,v)=\bigg(\frac{p}{p-1}\bigg)^p F-\frac{p^p}{p-1}\frac{f^p}{(A+(p-1)v)^{p-1}},
\end{equation*}
defined over the domain
\[\mathcal D:=\bigg \{ (F,f,A,v) \in \mathbb R^4 \mid F > 0, f > 0, A > 0, v > 0, v \geq A, f^p\leq F v^{p-1} \bigg\}.\]
We proved in section 1 that $\mathcal B$ satisfies the inequality (\ref{main inequality stronger}).
We are now going to show how the inequality (\ref{main inequality stronger}) entails a differential inequality for the function $\mathcal B$.\\
Let us consider a fixed point $(\tilde F, \tilde f, \tilde A, \tilde v)$ in the set of the interior points of $\mathcal D$. Let us consider $a\geq 0$, $b\geq 0$, $c\geq 0$. Let us consider $t>0$. We now define
\[\phi(t)=(\tilde F+(t b)^p,\tilde f+t^2 ab,\tilde A+t c,\tilde v+(t a)^q)\in \mathbb R^4,\]
\[\psi(t)=(\tilde F+u_1 t,\tilde f+u_2 t,\tilde A+u_3 t,\tilde v+u_4 t )\in \mathbb R^4.\]
As long as we choose $\tilde t \in \mathbb R^+$ small enough, we have that $\phi(t) \in \mathcal D$ and $\psi(t) \in \mathcal D$ for all  $0\leq t<\tilde t$.\\
So we may now compute the main inequality (\ref{main inequality stronger}) in the following way
\begin{align*}
&\mathcal B\bigg(\tilde F+(t^2 b)^p,\tilde f+(t^2)^2 ab,\tilde A+t^2 c,\tilde v+(t^2 a)^{p'}\bigg )-\mathcal B\bigg (\tilde F,\tilde f,\tilde A,\tilde v\bigg )+\\
&\mathcal B\bigg (\tilde F,\tilde f,\tilde A,\tilde v\bigg )-\frac{1}{2}\bigg[\mathcal B\bigg (\tilde F+u_1t,\tilde f+u_2t,\tilde A+u_3t,\tilde v+u_4t\bigg )+\\
&\mathcal B\bigg (\tilde F-u_1t,\tilde f-u_2t,\tilde A-u_3t,\tilde v-u_4t\bigg )\bigg]\geq p^p\frac{f^p}{(A+(p-1)v)^p}t^2 c,
\end{align*}
which is equivalent to
\begin{align*}
&\mathcal B(\phi(t^2) )-\mathcal B(\phi(0))+\mathcal B(\psi(0))-\frac{1}{2}\bigg[\mathcal B(\psi(-t))+\mathcal B(\psi(t))\bigg]\geq\\
& p^p\frac{(\tilde f+(t^2)^2 ab)^p}{(\tilde A+t^2 c+(p-1)(\tilde v+(t^2a)^{p'}))^p}t^2 c.
\end{align*}
We are allowed to compute this inequality because, by setting $X=\phi(t^2)$, $\tilde X=\phi(0)=\psi(0)$ and $X_+=\psi(t)$, $X_-=\psi(-t)$, we have 
\[X=\tilde X+\bigg((t^2b)^p,(t^2a)\cdot(t^2b),(ta)^{p'},c^2t\bigg),\]
\[\tilde X=\frac 1 2 (X_++X_-),\]
and $X$, $\tilde X$, $X_+$, $X_-$ are in the domain $\mathcal D$, so the hypotheses of the main inequality are satisfied.\\
Dividing by $t^2$ and taking the limit as $t \rightarrow 0$ we get
\[\lim_{t\rightarrow 0}\frac{\mathcal B(\phi(t^2) )-\mathcal B(\phi(0))+\mathcal B(\psi(0))-\frac{1}{2}\bigg[\mathcal B(\psi(-t))+\mathcal B(\psi(t))\bigg]}{t^2}\geq p^p\frac{\tilde f^p}{(\tilde A+(p-1)\tilde v)^p} c.\]
By a change of variable we get
\[\lim_{s\rightarrow 0}\frac{\mathcal B(\phi(s) )-\mathcal B(\phi(0))}{s}+\lim_{t\rightarrow 0} \frac{\mathcal B(\psi(0))-\frac{1}{2}\bigg[\mathcal B(\psi(-t))+\mathcal B(\psi(t))\bigg]}{t^2}\geq p^p\frac{\tilde f^p}{(\tilde A+(p-1)\tilde v)^p} c,\]
so we get
\begin{equation}
\frac{\partial}{\partial t}\mathcal B(\phi(t))\bigg|_{t=0}-\frac{1}{2}\frac{\partial ^2}{\partial t^2}\mathcal B(\psi(t))\bigg|_{t=0}\geq p^p\frac{\tilde f^p}{(\tilde A+(p-1)\tilde v)^p} c.
\end{equation}
By computing the derivative we get
\[\langle \nabla \mathcal B(\phi(0)),\phi'(0)\rangle -\frac{1}{2}\bigg[\langle \mathcal H (\mathcal B)(\psi(0)) \psi'(0),\psi'(0)\rangle+\langle \nabla \mathcal B(\psi(0)), \psi''(0)\rangle\bigg]\geq p^p\frac{\tilde f^p}{(\tilde A+(p-1)\tilde v)^p} c.\]
Now we observe that
\[\phi'(0)=(0,0,c,0), \quad \psi'(0)=(u_1,u_2,u_3,u_4)=:u, \quad \psi''(0)=(0,0,0,0).\]
So we get
\begin{equation}\label{continuous version of inequality}
\frac{\partial \mathcal B}{\partial x_3}\cdot  c -\frac{1}{2}\langle \mathcal H( \mathcal B)\cdot u, u\rangle \geq p^p\frac{\tilde f^p}{(\tilde A+(p-1)\tilde v)^p} c
\end{equation}
for any $c\geq0$.\\
We may verify that the function  
\begin{equation*}
\mathcal B(F,f,A,v)=\bigg(\frac{p}{p-1}\bigg)^p F-\frac{p^p}{p-1}\frac{f^p}{(A+v(p-1))^{p-1}}
\end{equation*}
satisfies the inequality (\ref{continuous version of inequality}). We compute
\[\frac{\partial \mathcal B}{\partial x_3}(\tilde F, \tilde f, \tilde A, \tilde v)\cdot c=p^p\frac{\tilde f^p}{(\tilde A + \tilde v(p-1))^p}c,\]
and $\mathcal B$ is concave so it satisfies $-\frac{1}{2}\langle \mathcal H( \mathcal B)\cdot u, u\rangle\geq 0$, showing that (\ref{continuous version of inequality}) is satisfied.\\
It follows that the function $\mathcal B$ satisfies the inequality
\[-\frac{\partial \mathcal B(x)}{\partial x_3}u_5 +\frac{1}{2}\sum_{i,j=1}^{4}\frac{\partial^2\mathcal B(x)}{\partial x_i \partial x_j}u_i u_j +p^p\bigg(\frac{x_2}{x_3+(p-1)x_4}\bigg)^pu_5 \leq0\quad \forall x \in \mathcal D, \; \forall u \in \mathbb R^5, \; u_5\geq0.\]
So the function $\mathcal B$ satisfies the following inequality
\begin{equation}\label{HJB supersolution}
\sup_{\underset{u_5\geq0}{u \in \mathbb R^5}}\bigg\{-\frac{\partial \mathcal B(x)}{\partial x_3}u_5 +\frac{1}{2}\sum_{i,j=1}^{4}\frac{\partial^2\mathcal B(x)}{\partial x_i \partial x_j}u_i u_j +p^p\bigg(\frac{x_2}{x_3+(p-1)x_4}\bigg)^pu_5\bigg\}\leq 0,
\end{equation}
so we will read the function $\mathcal B$ as a supersolution to a Hamilton-Jacobi-Bellman equation. Moreover, $\mathcal B$ is actually a solution to the Hamilton-Jacobi-Bellman equation by taking $u_1=u_2=u_3=u_4=0$.\\
So $\mathcal B$ satisfies the Hamilton-Jacobi-Bellman equation
\begin{equation}\label{HJB solution}
\sup_{\underset{u_5\geq0}{u \in \mathbb R^5}}\bigg\{-\frac{\partial \mathcal B(x)}{\partial x_3}u_5 +\frac{1}{2}\sum_{i,j=1}^{4}\frac{\partial^2\mathcal B(x)}{\partial x_i \partial x_j}u_i u_j +p^p\bigg(\frac{x_2}{x_3+(p-1)x_4}\bigg)^pu_5\bigg\}=0 \quad \forall x \in \mathcal D.
\end{equation}
So we naturally got a Hamilton-Jacobi-Bellman equation that can be interpreted as the equation associated to a stochastic optimal control problem.\\
\subsection{Stochastic optimal control problem}
Now we enunciate a stochastic optimal control which defines a Bellman function $g$ such that $g\equiv \mathcal B$.\\
Let us consider the following extension of the function $\mathcal B$ to the closure $\overline {\mathcal D}$ of its domain:
\[\tilde {\mathcal B}:\overline {\mathcal D} \longrightarrow \mathbb R\]
defined in the following way:
\[\tilde{\mathcal B}(x) = \begin{cases}
\mathcal B(x) \quad\quad\quad\quad\quad\,\, \text{if } x \in \mathcal D,\\
\underset{y \rightarrow x}{\liminf}\;\mathcal B(y) \quad \text{if } x \in \overline{\mathcal D} \backslash \mathcal D.\end{cases}\]
\begin{Rem}
For all points $x \in \overline{\mathcal D} \backslash \mathcal D$ such that $(x_3,x_4) \neq (0,0)$ the function $\mathcal B$ extends continuously to the value
\begin{equation*}
\tilde{\mathcal B}(x_1,x_2,x_3,x_4)=\underset{y\rightarrow x}{\lim}\mathcal B(y)=\bigg(\frac{p}{p-1}\bigg)^p x_1-\frac{p^p}{p-1}\frac{{x_2}^p}{(x_3+x_4(p-1))^{p-1}}.
\end{equation*}
The remaining points $x \in \overline{\mathcal D} \backslash \mathcal D$ are the points $x=(F,f,0,0)$, however by definition of $\mathcal D$ we have $f^p\leq F v^{p-1}$, so $f=0$. For the points $x=(F,0,0,0)$ such that $F\geq 0$ we have
\[\tilde{\mathcal B}(x)=\liminf_{y\rightarrow x} \mathcal B(y)=0. \]
\end{Rem}
\begin{proof}
We are going to show this fact by recalling that $\mathcal B \geq0$, so $\underset{y\rightarrow x}{\liminf}\; \mathcal B(y)\geq0$, and by considering a proper sequence of points. Let $v\geq A>0$, let $0<t\leq 1$. Let us first assume that $F>0$. Let us consider the points
\[x(t)=(F,(F(tv)^{p-1})^{\frac 1 p},t^2 A,tv).\]
By construction $x(t) \in \mathcal D$, $\underset{t\rightarrow 0}{\lim}\;x(t)=(F,0,0,0)$, and
\begin{align*}
\lim_{t\rightarrow 0} \mathcal B(x(t))=&\lim_{t\rightarrow 0}\bigg[\bigg(\frac{p}{p-1}\bigg)^p F-\frac{p^p}{p-1}\frac{F(tv)^{p-1}}{(t^2A+(p-1)tv)^{p-1}}\bigg]=\\
&F\cdot \lim_{t\rightarrow 0} \bigg[\bigg(\frac{p}{p-1}\bigg)^p -\frac{p^p}{p-1}\frac{v^{p-1}}{(tA+(p-1)v)^{p-1}}\bigg]=0.
\end{align*}
Let us assume $F=0$. Let  $\tilde F>0$. We consider the sequence
\[x(t)=(t\tilde F,(t\tilde F(tv)^{p-1})^{\frac 1 p},t^2 A,tv),\]
and the proof holds with the same argument.\\
So $0\leq\underset{y\rightarrow x}{\liminf}\; \mathcal B(y)\leq \underset{t\rightarrow 0}{\lim}\; \mathcal B(x(t))=0$, which ends the  proof.
\end{proof}

Let $x \in \mathcal D$, $t \geq0$, $u=(u_1,u_2,\dots,u_5) \in \mathbb R^5$ such that $u_5\geq0$. We will now define a payoff density and a bequest function to get the stochastic optimal control problem we are looking for. These functions will not depend on the time variable, so in the notation we will skip writing it. Let us define the payoff density
\[{\eta}^u(x,t)\equiv {\eta}^u(x):=p^p\bigg(\frac{x_2}{x_3+(p-1)x_4}\bigg)^p u_5.\]
Let $x \in \overline{\mathcal D}$. We define the bequest function
\[K(x_1,x_2,x_3,x_4,t)\equiv K(x):=\tilde {\mathcal B}(x),\]
i.e. $K$ is the function
\[ K(x_1,x_2,x_3,x_4,t)\equiv K(x)=\begin{cases}
\bigg(\frac{p}{p-1}\bigg)^px_1-\frac{p^p}{p-1}\frac{x_2^p}{(x_3+(p-1)x_4)^{p-1}} \quad \text{if } x \in \mathcal D,\\
\underset{y \rightarrow x}{\liminf}\;\mathcal B(y) \quad \quad \;\;\text{if } x \in \overline{\mathcal D} \backslash \mathcal D.\end{cases}\]
To finish the formulation of the stochastic optimal control problem we define the coefficients
\[b(u,x,t)\equiv b(u,x):=(0,0,-u_5,0),\]
\[\sigma (u,x,t)\equiv \sigma(u,x):=(u_1,u_2,u_3,u_4). \]
Let $\{u_t\}_{t\geq0}$ be a control such that $u_t(\omega) \in \{u \in \mathbb R^5 \; | \; u_5 \geq 0\}$. We consider the stochastic process $\{X_t\}=\{(F_t,f_t,A_t,v_t)\}$ solution to the stochastic differential equation (\ref{process bellman def}). The Bellman function associated to the problem is
\[g(x)=\sup_{\{u_t\}} E^x\bigg[\int_0^{\tau_{\mathcal D}} p^p\bigg(\frac{f_s}{A_s+(p-1)v_s}\bigg)^p u_5 \;ds +K(X_{\tau_{\mathcal D}})\chi_{\{\tau_{\mathcal D}<+\infty\}}\bigg],\]
where the supremum is taken over the set of controls $\{u_t\}_{t\geq0}$ such that $\{u_t\}$ is measurable with respect to $\mathcal F_t$, where $\{\mathcal F_t\}_{t\geq0}$  is the filtration generated by the variables $\{B_s \mid 0\leq s \leq t\}$, and such that the values $u_t(\omega)$ belong to the set $\{(u_1,u_2,u_3,u_4,u_5) \in \mathbb R^5 \mid u_5 \geq 0\}$.\\
So by the Hamilton-Bellman-Jacobi equation theorem (see \cite{Oksendal}, theorem 11.2.1) the function $g$ satisfies the equation (\ref{HJB solution}). We will also write the equation (\ref{HJB solution}) in the following way
\begin{equation}\label{HJB equation with generator}
\sup_{\underset{u_5\geq0}{u \in \mathbb R^5}}\bigg\{ (\mathcal L^ug)(x)+p^p\bigg(\frac{x_2}{x_3+(p-1)x_4}\bigg)^p u_5\bigg\}= 0 \quad \forall x \in \mathcal D.
\end{equation}
We recall that the operator 
\begin{equation}
\mathcal (\mathcal L^u\varphi)(x)=-\frac{\partial \varphi(x)}{\partial x_3}u_5 +\frac{1}{2}\sum_{i,j=1}^{4}\frac{\partial^2\varphi(x)}{\partial x_i \partial x_j}u_i u_j
\end{equation}
 is the infinitesimal generator of the process $\{X_t\}$ solution to the equation (\ref{process bellman def}) for the choice of the control $\{u_t\}$ such that $u_t\equiv u  \in \{y \in\mathbb R^5 \mid y_5\geq0\} $. Indeed, the infinitesimal generator $\mathcal A$ of such process (see \cite{Oksendal}, theorem 7.3.3) can be characterized by 
\[(\mathcal A g)(x)=\sum_{i=1}^{4}b_i(u,x)\frac{\partial g}{\partial x_i}(x)+\frac{1}{2}\sum_{i,j=1}^{4} (\sigma \sigma ^T)_{i,j}(u,x)\frac{\partial^2 g}{\partial x_i \partial x_j}(x)=(\mathcal L^ug)(x).\]

\subsection{The dyadic Bellman function is a stochastic Bellman function}
How we prove Theorem \ref{dyadic bellman is stochastic theorem}.
\begin{proof}
We are going to prove the stronger statement
\begin{equation}\label{tilde B is equal to stochastic bellman function}
g(x)=\tilde {\mathcal B}(x) \quad \forall x \in \overline{\mathcal D}.
\end{equation}
From  (\ref{tilde B is equal to stochastic bellman function}) and the definition of $\tilde {\mathcal B}$ it follows that
\begin{equation}
g(x)=\tilde {\mathcal B}(x)=\mathcal B(x) \quad \forall x \in \mathcal D,
\end{equation}
which is the required statement.
First we are going to prove that  $g(F,f,A,v)\geq\tilde{\mathcal B}(F,f,A,v)$. We compute
\[g(F,f,A,v) \geq E^{x_0}\bigg[ \int_0^{\tau_{\mathcal D}} {\eta}^{u_s}(X_s)ds+K(X_{\tau_{\mathcal D}})\chi_{\{\tau_{\mathcal D}<+\infty\}}\bigg]\]
for the choice
\[u_t=(0,0,0,0,1), \quad x_0=(F,f,A,v).\]
Let us first suppose $v> 0$. By computation we get $F_s \equiv F$, $f_s \equiv f$, $v_s \equiv v$, $A_s =A-s$, and $\tau_{\mathcal D}=A$, so, since the control is deterministic, we get
\begin{align*}
g(F,f,A,v)\geq& \int_0^{A} p^p\bigg(\frac{f}{A-s+(p-1)v}\bigg)^p \;ds +K(F,f,A-A,v)=\\
&\bigg[\frac{p^p}{p-1}\frac{f^p}{(A-s+(p-1)v)^{p-1}} \bigg]_{s=0}^{s=A}+\tilde {\mathcal B}(F,f,0,v)=\\
&\bigg(\frac{p}{p-1}\bigg)^p \frac{f^p}{v^{p-1}}-\frac{p^p}{p-1}\frac{f^p}{(A+(p-1)v)^{p-1}}+\liminf_{y\rightarrow (F,f,0,v)}\mathcal B(y)=\\
&\bigg(\frac{p}{p-1}\bigg)^p \frac{f^p}{v^{p-1}}-\frac{p^p}{p-1}\frac{f^p}{(A+(p-1)v)^{p-1}}+\bigg(\frac{p}{p-1}\bigg)^p F-\bigg(\frac{p}{p-1}\bigg)^p\frac{f^p}{v^{p-1}}=\\
&\bigg(\frac{p}{p-1}\bigg)^p F-\frac{p^p}{p-1}\frac{f^p}{(A+(p-1)v)^{p-1}}=\tilde{\mathcal B}(F,f,A,v).
\end{align*}
On the other hand, if $v=0$ then $A=0$ and $f=0$ by definition of the domain $\mathcal D$, so in this case the stopping time $\tau_{\mathcal D}$ is equal to 0, so the profit gain over the trajectory is 0, and we are left with the bequest gain. So the inequality becomes
\[g(F,f,A,v)\geq 0+K(F,0,0,0)=\tilde{\mathcal B}(F,0,0,0)=0,\]
which ends the proof that $g \geq \tilde{\mathcal B}$.\\

To prove that $g \leq \tilde{\mathcal B}$ we are going to first enunciate a heuristic argument to show it, using Jensen's inequality.\\
Let us consider a control $u=\{u_t\}_{t\geq0}$ such that $u_t=(u_{1}(t),\dots,u_4(t),0)$ for $0\leq t < s$ and then $u_t=(0,0,0,0,1)$ for $t\geq s$. Let $\{X_t\}$ be the solution to (\ref{process bellman def}) for this choice of the control $\{u_t\}$. The control $\{u_t\}$ lets the process $\{X_t\}$ behave like a martingale diffusion (the process has no drift) up to the time $s$, and on this part of the trajectory there is no profit gain (because the profit density is equal to 0 when $u_5=0$). Moreover, the control $\{u_t\}$ lets the process $\{X_t\}$ drift towards the boundary of the domain from the time $s$ onwards.\\
Let $t \mapsto X(\omega)(t):=X_t(\omega)$ be a trajectory of the process $\{X_t\}$.\\
If $\tau_D(\omega)\leq s$ then, by continuity of the process $\{X_t\}$, the trajectory $X(\omega)$ lands on the point $X_{\tau_{\mathcal D}(\omega)}(\omega) \in \partial \mathcal D$ for almost all the $\omega$ with such properties. So almost all trajectories $X(\omega)$ such that $\tau_{\mathcal D}(\omega)\leq s$ gain an amount of profit equal to $K(X_{\tau_{\mathcal D}(\omega)})=\tilde{\mathcal B}(X_{\tau_{\mathcal D}(\omega)})=\tilde{\mathcal B}(X_{s\wedge \tau_{\mathcal D}(\omega)})$.\\
If  $\tau_D(\omega)>s$, the trajectory $t \mapsto X_t(\omega)$ of the process $\{X_t\}$ lands on a point $x$ in the interior of the domain $\mathcal D$ at the time $s$, without exiting the domain $\mathcal D$ before the time $s$. We observe that, in the first part of the proof, we proved that a control $\{\hat {u}_t\}$ such that $\hat{u}_t=(0,0,0,0,1)$ generates a process $\{\hat{X}_t\}$ that gains an amount of average profit equal to the value of the function $\tilde{\mathcal B}$ in the starting point. So, by this observation, it follows that the trajectory $X(\omega)$ gains no profit during the time $0<t<s$, and gains an amount of profit equal to $\tilde{\mathcal B}(X_{s})=\tilde{\mathcal B}(X_{\tau_{\mathcal D}(\omega) \wedge s})$ during the times $s\leq t \leq \tau_{\mathcal D}(\omega)$.\\
Based on these observations the average amount of profit gained by $\{X_t\}$, given a control $\{u_t\}$ of this kind, is
\[J^u(x)=E^x\bigg[0+\tilde{\mathcal B}(X_{s\wedge \tau_{\mathcal D}})\bigg],\]
where the first addent stands for the null gain on the trajectory up to the time $s\wedge \tau_{\mathcal D}$, while the second addend is equal to the gain from that moment onwards (the profit gain from the bequest function at the end of times is included in the second addend). \\
Following this notation the Bellman function $g$ is
\[g(x)=\sup_{u=\{u_t\}}J^u(x).\]
However, $\tilde{\mathcal B}$ is concave, so by Jensen's inequality we have
\[E^x\bigg[\tilde{\mathcal B}(X_{s\wedge \tau_{\mathcal D}})\bigg]\leq \tilde{\mathcal B}\bigg(E^x[X_{s\wedge \tau_{\mathcal D}}]\bigg).\]
Moreover, $\{X_t\}$ is a martingale up to the time $s$ by definition, so we get
\[J^u(x)=E^x\bigg[\tilde{\mathcal B}(X_{s\wedge \tau_{\mathcal D}})\bigg]\leq  \tilde{\mathcal B}\bigg(E^x[X_{s\wedge \tau_{\mathcal D}}]\bigg)= \tilde{\mathcal B}(x),\]
and the heuristic idea is that there is ``independence" between letting the process drift in the third variable (which gives a non-negative gain) and letting the process be a diffusion (which gains nothing), so we can let the process be a combination of the two and the argument will still hold. So by taking the supremum over all controls $\{u_t\}$ we get
\[g(x)=\sup_{u=\{u_t\}}J^u(x)\leq \tilde{\mathcal B}(x).\]

We are now going to give a proof that $g \leq \tilde{\mathcal B}$ using Dynkin's formula (see \cite{Oksendal}, theorem 7.4.1).\\
We will skip some technical details in the following proof.\\
Let $\{u_t\}_{t\geq0}$ be a given control. Let $\{X_t\}_{t\geq0}$ be the process solution to (\ref{process bellman def}) for this choice of the control $\{u_t\}_{t\geq0}$. Let $\tau_{\mathcal D}$ be the first exit time for $\{X_t\}$ from $\mathcal D$. We will first assume that $\tau_{\mathcal D}<+\infty$ almost surely. Now we apply Dynkin's formula
\begin{equation}\label{Dynkin's formula}
E^x[\tilde{\mathcal B}(X_{\tau_{\mathcal D}})]=\tilde{\mathcal B}(x) + E^x \bigg[ \int_0^{\tau_{\mathcal D}} (\mathcal L^{u_s} \tilde{\mathcal B})(X_s)dx\bigg]
\end{equation} 
to the function $\tilde{\mathcal B}$ and the process $\{X_t\}_{t\geq 0}$. We get
\[\tilde{\mathcal B}(x)=E^x[\tilde{\mathcal B}(X_{\tau_{\mathcal D}})] - E^x \int_0^{\tau_{\mathcal D}} (\mathcal L^{u_s}\tilde{\mathcal B})(X_s)ds.\]
Now, since $\tau_{\mathcal D}<+\infty$ almost surely, the event $\chi_{\{\tau_{\mathcal D}<+\infty\}}$ has a probability of 1, so
\[\tilde{\mathcal B}(x)=E^x[\tilde{\mathcal B}(X_{\tau_{\mathcal D}})\chi_{\{\tau_{\mathcal D}<+\infty\}}] - E^x \int_0^{\tau_{\mathcal D}} (\mathcal L^{u_s}\tilde{\mathcal B})(X_s)ds.\]
The equation (\ref{HJB equation with generator}) entails that $-(\mathcal L^{u_s}\tilde{\mathcal B})(y)\geq {\eta}^{u_s}(y)$, so we get
\[\tilde{\mathcal B}(x)\geq E^x \int_0^{\tau_{\mathcal D}} {\eta}^{u_s}(X_s)ds +E^x[\tilde{\mathcal B}(X_{\tau_{\mathcal D}})\chi_{\{\tau_{\mathcal D}<+\infty\}}].\]
However, $\tilde{\mathcal B}(X_{\tau_{\mathcal D}})=K(X_{\tau_{\mathcal D}})$, so we get
\[\tilde{\mathcal B}(x)\geq E^x \int_0^{\tau_{\mathcal D}} {\eta}^{u_s}(X_s)ds +E^x[K(X_{\tau_{\mathcal D}})\chi_{\{\tau_{\mathcal D}<+\infty\}}]=g(x).\]
If $\tau_{\mathcal D}$ is not almost surely finite, we are going to show an idea of the proof. We may consider the stopping time $\tau(T)=\tau_{\mathcal D}\wedge T=\min\{\tau_{\mathcal D},T\}$ for $T>0$. This procedure is equivalent to considering the processes $Y_t=(t,X_t)$ in the domain $[0,T]\times \mathcal D$, and then defining the Bellman function $\tilde{\mathcal B}_T$ associated to those processes, which is a standard way to define the Bellman functions (see \cite{Oksendal}, chapter 11). \\
Since $\tau(T)<+\infty$ almost surely, we may apply Dynkin's formula to that stopping time and, with the same argument we used before, we get
\begin{align*}
\tilde{\mathcal B}_T(&x)\geq E^x \int_0^{\tau(T)} {\eta}^{u_s}(X_s)ds +E^x[K(X_{\tau(T)})\chi_{\{\tau(T)<+\infty\}}]\\
&\downarrow T \rightarrow + \infty \quad\quad\quad\quad \downarrow  T \rightarrow + \infty\\
\tilde{\mathcal B}(&x) \geq E^x \int_0^{\tau_{\mathcal D}} {\eta}^{u_s}(X_s)ds +E^x[K(X_{\tau_{\mathcal D}})\chi_{\{\tau_{\mathcal D}<+\infty\}}].
\end{align*}
So, by taking the supremum over all controls $\{u_t\}$, we get
\[\tilde{\mathcal B}(x) \geq \sup_{\{u_t\}}E^x\bigg[ \int_0^{\tau_{\mathcal D}} {\eta}^{u_s}(X_s)ds+K(X_{\tau_{\mathcal D}})\chi_{\{\tau_{\mathcal D}<+\infty\}}\bigg]=g(x),\]
which ends the proof that $\tilde{\mathcal B}\equiv v$, so $\tilde{\mathcal B}$ is the Bellman function solution to the stochastic optimal control problem.
\end{proof}

\section{Appendix}
We are going to prove the following lemma.
\begin{Lem}
The domain
\[\mathcal D:=\bigg \{ (F,f,A,v) \in \mathbb R^4 \mid F>0, f >0, A > 0, v > 0, v \geq A, f^p\leq F v^{p-1} \bigg\}\]
is convex. 
\end{Lem}
\begin{proof}
We write the domain $\mathcal D$ in the form
\[\mathcal D=\{v\geq A\}\cap \mathcal A, \]
here $\mathcal A$ is the set
\[\mathcal A=\{(F,f,A,v)\in\mathbb R^4 \mid F>0, f>0, v >0, f^p\leq Fv^{p-1}\}.\]
To prove that the domain is convex we just need to prove that it is an intersection of convex sets.\\ 
The set $\{v\geq A\}$ is trivially convex because it is a half-plane.
Since $\frac{1}{p}+\frac{1}{p'}=1$ and $\frac{p-1}{p}=\frac{1}{p'}$, the set $\mathcal A$ can be written in the form
\[\mathcal A=\mathcal S \cap \{f > 0\},\]
here $\{f >0\}$ is another  half-plane (a convex set), while $\mathcal S$ is the set
\[\mathcal S=\{(F,f,A,v)\in\mathbb R^4 \mid F> 0, v >0,  f\leq F^{\frac{1}{p}}v^{\frac{1}{p'}} \}.\]
The set $\mathcal S$ is the subgraph of the function
\begin{align*}
h:\mathbb R^+ \times \mathbb R \times \mathbb R^+ &\longrightarrow \mathbb R_0^+\\
(F,A,v) &\longmapsto   F^{\frac{1}{p}}v^{\frac{1}{p'}}.
\end{align*}
To prove that $\mathcal S$ is convex, all we need to do is to prove that $h$ is a concave function (since $h$ is defined over a convex domain).\\
Since $h$ does not depend on the variable $A$, we will treat it as a function over the other two variables only:
\begin{align*}
h:\mathbb R^+  \times \mathbb R^+ &\longrightarrow \mathbb R^+\\
(F,v) &\longmapsto   F^{\frac{1}{p}}v^{\frac{1}{p'}}.
\end{align*}
We compute the Hessian matrix of the function $h$: for all $F>0$, $v>0$
\begin{align*}
\frac{\partial h}{\partial F}(F,v)=\frac{1}{p}F^{\frac{1}{p}-1}v^{\frac{1}{p'}}, &\quad \frac{\partial h}{\partial v}(F,v)=\frac{1}{p'}F^{\frac{1}{p}}v^{\frac{1}{p'}-1}.\\
\\
\frac{\partial^2 h}{\partial F^2}(F,v)=\frac{1-p}{p^2}F^{\frac{1}{p}-2}v^{\frac{1}{p'}}, &\quad \frac{\partial^2 h}{\partial v\partial F}(F,v)=\frac{1}{p{p'}}F^{\frac{1}{p}-1}v^{\frac{1}{p'}-1},\\
 \frac{\partial^2 h}{\partial F \partial v}(F,v)\frac{1}{p{p'}}F^{\frac{1}{p}-1}v^{\frac{1}{p'}-1},& \quad \frac{\partial^2 h}{\partial v^2}(F,v)=\frac{1-{p'}}{{p'}^2}F^{\frac{1}{p}}v^{\frac{1}{p'}-2}.\\
\end{align*}
So the Hessian matrix is
\begin{equation}\label{hessian 1}
\mathcal H (h)(F,v)=\begin{bmatrix}
\frac{1-p}{p^2}F^{\frac{1}{p}-2}v^{\frac{1}{p'}}&\frac{1}{p{p'}}F^{\frac{1}{p}-1}v^{\frac{1}{p'}-1}\\
\\
 \frac{1}{p{p'}}F^{\frac{1}{p}-1}v^{\frac{1}{p'}-1}&\frac{1-{p'}}{{p'}^2}F^{\frac{1}{p}}v^{\frac{1}{p'}-2}
\end{bmatrix}.
\end{equation}
If the Hessian matrix of $h$ has non-positive eigenvalues then the function $h$ is concave.\\
Now we compute the eigenvalues of the Hessian matrix (\ref{hessian 1}):
\begin{align*}
\det (\mathcal H (h)(F,v)-\lambda I)=&\det \begin{bmatrix}
\frac{1-p}{p^2}F^{\frac{1}{p}-2}v^{\frac{1}{p'}}-\lambda&\frac{1}{p{p'}}F^{\frac{1}{p}-1}v^{\frac{1}{{p'}}-1}\\
\\
 \frac{1}{p{p'}}F^{\frac{1}{p}-1}v^{\frac{1}{p'}-1}&\frac{1-{p'}}{{p'}^2}F^{\frac{1}{p}}v^{\frac{1}{p'}-2}-\lambda
\end{bmatrix}=\\
&\frac{(1-p)(1-{p'})}{(p{p'})^2}F^{\frac{1}{p}+\frac{1}{p}-2}v^{\frac{1}{p'}+\frac{1}{p'}-2}-\frac{1}{(p{p'})^2}F^{2(\frac{1}{p}-1)}v^{2(\frac{1}{p'}-1)}-\\
&\lambda \bigg[\frac{1-p}{p^2}F^{\frac{1}{p}-2}v^{\frac{1}{p'}}+\frac{1-{p'}}{{p'}^2}F^{\frac{1}{p}}v^{\frac{1}{p'}-2}\bigg]+ \lambda^2.
\end{align*}
Now we recall that $p{p'}=p+{p'}$, so $(1-p)(1-{p'})=1-p-{p'}+p{p'}=1-p-{p'}+p+{p'}=1$, so we get
\begin{align*}
\det (\mathcal H (h)(F,v)-\lambda I)=\lambda^2-\lambda \bigg[\frac{1-p}{p^2}F^{\frac{1}{p}-2}v^{\frac{1}{p'}}+\frac{1-{p'}}{{p'}^2}F^{\frac{1}{p}}v^{\frac{1}{p'}-2}\bigg].
\end{align*}
The eigenvalues of $\mathcal H (h)(F,v)$ are the solutions to the following equation equation of variable $\lambda$:
\[\det (\mathcal H (h)(F,v)-\lambda I)=0.\]
The solutions are the two values 
\[\lambda_1=0,\quad \lambda_2=\frac{1-p}{p^2}F^{\frac{1}{p}-2}v^{\frac{1}{p'}}+\frac{1-{p'}}{{p'}^2}F^{\frac{1}{p}}v^{\frac{1}{p'}-2}.\]
Now we observe that $1-p<0$, $1-{p'}<0$ and $F>0$, $v>0$, so the second eigenvalue is $\lambda_2<0$, so the Hessian matrix $\mathcal H (h)(F,v)$ is negative semi-definite for all $F>0$ and $v>0$, so this entails that $h$ is concave, and the subgraph $\mathcal S$ is a convex set. So the domain $\mathcal D$\ of the function $\mathcal B$ in (\ref{bellman function}) is a convex set since it is an intersection of convex sets.
\end{proof}

We are going to prove the following lemma.
\begin{Lem}
The function
\begin{equation*}
\mathcal B(F,f,A,v)=\bigg(\frac{p}{p-1}\bigg)^p F-\frac{p^p}{p-1}\frac{f^p}{(A+(p-1)v)^{p-1}}
\end{equation*}
is concave. 
\end{Lem}
\begin{proof}
We will compute the eigenvalues of the Hessian Matrix $\mathcal H (\mathcal B)$. We will compute the actual $4 \times 4$ Hessian matrix (without reducing it to a $2 \times 2$ matrix), because the computation can be useful to compute the eigenvectors, which may be useful to study the properties of some of the stochastic processes associated to the problem.\\
In the following equations we omit writing the dependence of the derivatives from the variables $(F,f,A,v)$ to simplify the notations. The first order derivatives are
\begin{align*}
&\frac{\partial \mathcal B}{\partial F}=\bigg(\frac{p}{p-1}\bigg)^p, &\frac{\partial \mathcal B}{\partial f}=-\frac{p^{p+1}}{p-1} \frac{f^{p-1}}{(A+(p-1)v)^{p-1}},\\
&\frac{\partial \mathcal B}{\partial A}=p^p\frac{f^{p}}{(A+(p-1)v)^{p}}, &\frac{\partial \mathcal B}{\partial v}=(p-1)p^p\frac{f^{p}}{(A+(p-1)v)^{p}}.
\end{align*}
The second order derivatives make up the rows of the Hessian matrix $\mathcal H(\mathcal B)$.\\
The first row is
\begin{align*}
&\frac{\partial^2 \mathcal B}{\partial F^2}=0, &\frac{\partial^2 \mathcal B}{\partial F\partial f}=0,\\
&\frac{\partial^2 \mathcal B}{\partial F\partial A}=0, &\frac{\partial^2 \mathcal B}{\partial F\partial V}=0.
\end{align*}
The second row is
\begin{align*}
&\frac{\partial^2 \mathcal B}{\partial f \partial F}=0, &\frac{\partial^2 \mathcal B}{\partial f^2}=-p^{p+1}\frac{f^{p-2}}{(A+(p-1)v)^{p-1}},	\\
&\frac{\partial^2 \mathcal B}{\partial f\partial A}=p^{p+1}\frac{f^{p-1}}{(A+(p-1)v)^{p}}, &\frac{\partial^2 \mathcal B}{\partial f\partial V}=(p-1)p^{p+1}\frac{f^{p-1}}{(A+(p-1)v)^{p}}.
\end{align*}
The third row is
\begin{align*}
&\frac{\partial^2 \mathcal B}{\partial A \partial F}=0, &\frac{\partial^2 \mathcal B}{\partial A \partial f}=p^{p+1}\frac{f^{p-1}}{(A+(p-1)v)^{p}},\\
&\frac{\partial^2 \mathcal B}{\partial A^2}=-p^{p+1}\frac{f^{p}}{(A+(p-1)v)^{p+1}}, &\frac{\partial^2 \mathcal B}{\partial A \partial V}=-(p-1)p^{p+1}\frac{f^{p}}{(A+(p-1)v)^{p+1}}.
\end{align*}
The fourth row is
\begin{align*}
&\frac{\partial^2 \mathcal B}{\partial V\partial F}=0, &\frac{\partial^2 \mathcal B}{\partial V \partial f}=(p-1)p^{p+1}\frac{f^{p-1}}{(A+(p-1)v)^{p}},\\
&\frac{\partial^2 \mathcal B}{\partial V \partial A}=-(p-1)p^{p+1}\frac{f^{p}}{(A+(p-1)v)^{p+1}}, &\frac{\partial^2 \mathcal B}{\partial V^2}=-(p-1)^2 p^{p+1}\frac{f^{p}}{(A+(p-1)v)^{p+1}}.
\end{align*}
\\
So the Hessian matrix of $\mathcal B$ at a point $(F,f,A,v)$ is 
\[\mathcal H (\mathcal B)=\begin{pmatrix}
0 & 0 & 0 & 0\\
0&-p^{p+1}\frac{f^{p-2}}{(A+(p-1)v)^{p-1}}&p^{p+1}\frac{f^{p-1}}{(A+(p-1)v)^{p}}&(p-1)p^{p+1}\frac{f^{p-1}}{(A+(p-1)v)^{p}}\\
0&p^{p+1}\frac{f^{p-1}}{(A+(p-1)v)^{p}}&-p^{p+1}\frac{f^{p}}{(A+(p-1)v)^{p+1}}&-(p-1)p^{p+1}\frac{f^{p}}{(A+(p-1)v)^{p+1}}\\
0 &(p-1)p^{p+1}\frac{f^{p-1}}{(A+(p-1)v)^{p}}&-(p-1)p^{p+1}\frac{f^{p}}{(A+(p-1)v)^{p+1}}&-(p-1)^2 p^{p+1}\frac{f^{p}}{(A+(p-1)v)^{p+1}}
\end{pmatrix}.\]
Let us compute the eigenvalues:
\begin{align*}
&0=\det(\mathcal H(\mathcal B)-\lambda I_4)=\\
\ \\
&-\lambda \det \begin{pmatrix}
-p^{p+1}\frac{f^{p-2}}{(A+(p-1)v)^{p-1}}-\lambda &p^{p+1}\frac{f^{p-1}}{(A+(p-1)v)^{p}}&(p-1)p^{p+1}\frac{f^{p-1}}{(A+(p-1)v)^{p}}\\
p^{p+1}\frac{f^{p-1}}{(A+(p-1)v)^{p}}&-p^{p+1}\frac{f^{p}}{(A+(p-1)v)^{p+1}}-\lambda &-(p-1)p^{p+1}\frac{f^{p}}{(A+(p-1)v)^{p+1}}\\
0&(p-1)\lambda &-\lambda
\end{pmatrix}=\\
\ \\
&\lambda^2 \bigg[ (p-1)\det \begin{pmatrix}
-p^{p+1}\frac{f^{p-2}}{(A+(p-1)v)^{p-1}}-\lambda&(p-1)p^{p+1}\frac{f^{p-1}}{(A+(p-1)v)^{p}}\\
p^{p+1}\frac{f^{p-1}}{(A+(p-1)v)^{p}} &-(p-1)p^{p+1}\frac{f^{p}}{(A+(p-1)v)^{p+1}}
\end{pmatrix}+\\
&\det \begin{pmatrix}
-p^{p+1}\frac{f^{p-2}}{(A+(p-1)v)^{p-1}}-\lambda&p^{p+1}\frac{f^{p-1}}{(A+(p-1)v)^{p}} \\
p^{p+1}\frac{f^{p-1}}{(A+(p-1)v)^{p}}&-p^{p+1}\frac{f^{p}}{(A+(p-1)v)^{p+1}}-\lambda
\end{pmatrix}\bigg]=\\
\ \\
&\lambda^2 \bigg[\lambda^2+p^{2p+2}\bigg[\frac{f^{p-2}}{(A+(p-1)v)^{p-1}}+\frac{f^{p}}{(A+(p-1)v)^{p+1}}\bigg]\lambda+\\
&p^{2p+2}(p-1)\frac{f^{p}}{(A+(p-1)v)^{p+1}} \lambda\bigg]=\\
&\lambda^3\bigg[\lambda +p^{2p+2}\bigg[\frac{f^{p-2}}{(A+(p-1)v)^{p-1}}+p\frac{f^{p}}{(A+(p-1)v)^{p+1}}\bigg]\bigg].
\end{align*}
So the eigenvalues of the Hessian matrix $\mathcal H (\mathcal B)$ are 0 of algebraic multiplicity 3 and \\$\tilde \lambda=-p^{2p+2}\bigg[\frac{f^{p-2}}{(A+(p-1)v)^{p-1}}+p\frac{f^{p}}{(A+(p-1)v)^{p+1}}\bigg]$ of algebraic multiplicity 1. However, since $f>0$, $v>0$ and $A>0$, then $\tilde \lambda > 0$, so all the eigenvalues are lower than or equal to 0. This entails that the Hessian matrix is negative semi-definite, so the function $\mathcal B$ is concave.
\end{proof}

We are going to prove the following lemma.
\begin{Lem}
The function
\[\mathcal B(F,f,A,v)=\bigg(\frac{p}{p-1}\bigg)^p F-\frac{p^p}{p-1}\frac{f^p}{(A+(p-1)v)^{p-1}}\]
satisfies \[0\leq\mathcal B(F,f,A,v) \leq \big(p/(p-1)\big)^p F\] for all $(F,f,A,v) \in \mathcal D$.
\end{Lem}
\begin{proof} 
The statement follows from the definition of the domain of the function. Since $p>1$, $F>0$, $f > 0$, $A>0$, $ v \geq A$ and $ f^p\leq F v^{p-1}$, we get
\begin{align*}
\mathcal B(F,f,A,v)=&\bigg(\frac{p}{p-1}\bigg)^p F-\frac{p^p}{p-1}\frac{f^p}{(A+(p-1)v)^{p-1}} \geq\\
& \bigg(\frac{p}{p-1}\bigg)^p F-\frac{p^p}{p-1}\frac{f^p}{(0+(p-1)v)^{p-1}}\geq \\
& \bigg(\frac{p}{p-1}\bigg)^p F-\bigg(\frac{p}{p-1}\bigg)^p\frac{F v^{p-1}}{v^{p-1}}\geq\\
& \bigg(\frac{p}{p-1}\bigg)^p\bigg[F-F\bigg]= 0,
\end{align*}
and
\[\mathcal B(F,f,A,v)=\bigg(\frac{p}{p-1}\bigg)^p F-\frac{p^p}{p-1}\frac{f^p}{(A+(p-1)v)^{p-1}}\leq \bigg(\frac{p}{p-1}\bigg)^p F.\]
\end{proof}
\section*{Acknowledgements}
\addcontentsline{toc}{section}{Acknowledgements}
We would like to thank Oliver Dragi\v cevi\'c and Nikolaos Chalmoukis for our useful talks about the topics in this article. The topics in this article constitute a part of the autor's PhD thesis ``Potential theory on metric spaces", which was written under the supervision of Nicola Arcozzi.\\

\end{document}